\documentclass{amsart}
\usepackage{amsmath,amsthm,amscd,amssymb}
\usepackage{latexsym}
\usepackage[dvips]{epsfig}
\allowdisplaybreaks
\newcommand{\Hmm}[1]{\leavevmode{\marginpar{\tiny%
$\hbox to 0mm{\hspace*{-0.5mm}$\leftarrow$\hss}%
\vcenter{\vrule depth 0.1mm height 0.1mm width \the\marginparwidth}%
\hbox to 0mm{\hss$\rightarrow$\hspace*{-0.5mm}}$\\\relax\raggedright #1}}}

\def\sign{\rm sign}
\def\Laplace{\Delta}

\def\bm{\left[ \begin{array}{cc}}
\def\endm{\end{array}\right]}
\def\les{\le\, C}
\def\beeq{\begin{equation}}
\def\eneq{\end{equation}}
\def\bea{\begin{eqnarray}}
\def\eea{\end{eqnarray}}

\def\R{\mathbb R}
\def\C{\mathbb C}

\def\spec{\rm spec}

\def\half{\frac12}
\def\eps{\epsilon}

\def\calB{\mathcal B}

\def\nn{\nonumber}
\def\Hil{{\mathcal H}}
\def\Dom{{\rm Dom}}
\def\la{\langle}
\def\ra{\rangle}
\def\Ran{{\rm Ran}}

\newcommand{\fint}{{\int}{\kern-1em \diagup}}
\newcommand{\pref}[1]{{\rm(\ref{#1})}}

\newcommand{\calH}{{\mathcal H}}
\newcommand{\calS}{{\mathcal S}}
\newcommand{\calF}{{\mathcal F}}
\newcommand{\calG}{{\mathcal G}}
\newcommand{\mr}{{\mathcal{R}}}

\newcommand{\f}{\frac}

\newcommand{\rank}{\text{\rm{rank}}}
\newcommand{\ran}{\text{\rm{Ran}}}

\newcommand{\beq}{\begin{equation}}
\newcommand{\eeq}{\end{equation}}
\newcommand{\ba}{\begin{align}}
\newcommand{\ea}{\end{align}}

\newcommand{\lang}{\left\langle}
\newcommand{\rang}{\right\rangle}

\DeclareMathOperator{\Ima}{Im}

\newtheorem{theorem}{Theorem}
\newtheorem*{theorem*}{Theorem}
\newtheorem{proposition}[theorem]{Proposition}
\newtheorem{lemma}[theorem]{Lemma}
\newtheorem*{lemma*}{Lemma}
\newtheorem{cor}[theorem]{Corollary}
\newtheorem*{corollary*}{Corollary}
\theoremstyle{definition}

\newtheorem*{definition*}{Definition}

\theoremstyle{remark}
\newtheorem{remark}{Remark}

\begin{document}

\title[Dispersive estimates]{Dispersive estimates for Schr\"{o}dinger operators in the presence of a resonance and/or
an eigenvalue at zero energy in dimension three: II}

\author{M. Burak Erdo\smash{\u{g}}an and Wilhelm Schlag}
\thanks{This work was initiated in June of 2004,
while the first author visited Caltech and he  wishes to thank
that institution for its hospitality and support. The first author
was partially supported by the NSF grant DMS-0303413. The second
author was partially supported by a Sloan fellowship and the NSF
grant DMS--0300081. The authors thank Avy Soffer for his interest
in this work.}

\address{Department of Mathematics \\
University of Illinois \\
Urbana, IL 61801, U.S.A.}

\address{Department of Mathematics \\
The University of Chicago \\
5734 South University Avenue\\
Chicago, IL 60637, U.S.A.}

\date{\today}

\begin{abstract}
We investigate boundedness of the evolution $e^{it\Hil}$ in the
sense of $L^2(\R^3)\to L^2(\R^3)$ as well as $L^1(\R^3)\to
L^\infty(\R^3)$ for the non-selfadjoint operator
\[ \Hil =
\left[\begin{array}{cc} -\Delta + \mu-V_1 & -V_2\\
V_2 & \Delta - \mu + V_1
\end{array}
\right]
\]
where $\mu>0$ and $V_1,V_2$ are real-valued decaying potentials.
Such operators arise when linearizing a focusing NLS equation around
a standing wave and the aforementioned bounds are needed in the
study of nonlinear asymptotic stability of such standing waves. We
derive our results under some natural spectral assumptions
(corresponding to a ground state soliton of NLS), see A1)--A4)
below, but without imposing any restrictions on the edges $\pm\mu$
of the essential spectrum. Our goal is to develop an ``axiomatic
approach", which frees the linear theory from any nonlinear context
in which it may have arisen.
\end{abstract}

\maketitle

\section{The matrix case: Introduction}
\label{sec:intro}

Consider the Schr\"{o}dinger operator  $H=-\Delta+V$ in $\R^3$,
where $V$ is a real-valued potential. Let $P_{ac}$ be the
orthogonal projection onto the absolutely continuous subspace of
$L^2(\R^3)$ which is determined by $H$. In Journ\'e, Soffer,
Sogge~\cite{JSS}, Yajima~\cite{Y1}, Rodnianski,
Schlag~\cite{RodSch}, Goldberg, Schlag~\cite{goldbergschlag} and
Goldberg~\cite{goldberg}, $L^1(\R^3)\rightarrow L^\infty(\R^3)$
dispersive estimates for the time evolution $e^{itH}P_{ac}$ were
investigated under various decay assumptions on the potential $V$
and the assumption that zero is neither an eigenvalue nor a
resonance of~$H$. Recall that zero energy is a resonance iff there
is $f\in L^{2,-\sigma}(\R^3)\setminus L^2(\R^3)$ for all
$\sigma>\half$ so that $Hf=0$. Here $L^{2,-\sigma}=\langle
x\rangle^{\sigma}L^2$ are the usual weighted $L^2$ spaces and
$\langle x\rangle:=(1+|x|^2)^{\frac12}$. For a  survey of recent
work in this area see~\cite{Sch2}.

In \cite{ErdSch}, the authors investigated dispersive estimates
when there is a resonance or eigenvalue at energy zero. It is
well-known, see Rauch~\cite{Rauch}, Jensen, Kato~\cite{JenKat},
and Murata~\cite{Mur}, that the decay in that case
is~$t^{-\half}$. Moreover, these authors derived expansions of the
evolution into inverse powers of time in weighted $L^2(\R^3)$
spaces. In \cite{ErdSch}, the authors obtained such expansions
with respect to the $L^1\to L^\infty$ norm, albeit only in terms
of the powers $t^{-\half}$ and $t^{-\frac32}$. Independently,
Yajima~\cite{Y2} achieved similar results.

In this paper we obtain analogous expansions for a class of {\bf
matrix Schr\"{o}dinger operators}. Consider the matrix
Schr\"{o}dinger operator
$$
\calH = \calH_0 + V =
\left[
\begin{array}{cc}
-\Delta + \mu & 0 \\
0 & \Delta - \mu
\end{array}\right] +
\left[\begin{array}{cc} -V_1 & -V_2\\
V_2 & V_1
\end{array}
\right]
$$
on $L^2(\R^3)\times L^2(\R^3)$. Here $\mu>0$ and  $V_1$, $V_2$ are
real-valued. It follows from Weyl's criterion that the essential
spectrum of $\calH$ is $(-\infty,-\mu]\cup[\mu,\infty)$. The
discrete spectrum may intersect $\C\setminus\R$, and the algebraic
and geometric multiplicities of eigenvalues may be different
(i.e., $\calH$ has a nonzero nilpotent part at these eigenvalues).

Such operators appear naturally as linearizations of a nonlinear
Schr\"odinger equation around a standing wave (or {\em soliton}),
see below. Dispersive estimates in the context of such
linearizations were obtained in Cuccagna~\cite{Cuc}, Rodnianski,
Schlag, Soffer~\cite{RSS1}, and \cite{Sch1} under various decay
assumptions on the potential and the assumption that zero is neither
an eigenvalue nor a resonance of~$\calH$. In addition, one always
assumes that there are no imbedded eigenvalues in the essential
spectrum.

The emphasis of the present paper is to develop an "abstract" (or
"axiomatic") approach, which frees the linear theory from any
reference to a nonlinear context in which it may have arisen. More
specifically, our results will require the following assumptions on
$\calH$ (in what follows,
$\sigma_3$ is one of the Pauli matrices, see~\eqref{eq:pauli}):\\

{\bf Assumptions}: \\
A1) $-\sigma_3 V$ is a positive matrix\\
A2) $L_-:=-\Delta + \mu - V_1 + V_2 \geq 0$ \\
A3) For some $\beta>0$,
\bea \label{decay}
|V_1(x)|+|V_2(x)|\lesssim \langle x\rangle^{-\beta} \eea
A4) There
are no imbedded eigenvalues in $(-\infty,-\mu)\cup(\mu,\infty)$

\bigskip

Assumptions A1)-A3) hold in the important example of a linearized
nonlinear Schr\"odinger equation, provided the linearization is
performed around the (positive) ground state standing wave. Indeed,
suppose that $\psi(t,x)=e^{it\alpha^2}\phi(x)$ is a standing wave
solution of the NLS
\[
i\partial_t\psi + \Delta\psi + |\psi|^{2\beta}\psi =0,
\]
where $\beta>0$.  Here we assume that $\phi$ is a ground state,
i.e.,
\[
\alpha^2 \phi - \Delta\phi = \phi^{2\beta+1}, \qquad \phi>0.
\]
Is known that such $\phi$ exist and that they are radial, smooth,
and exponentially decaying, see Strauss~\cite{Str1}, Berestycki,
Lions~\cite{BL1} and for uniqueness, see Coffman~\cite{Cof}, McLeod,
Serrin~\cite{McLSer}, and Kwong~\cite{Kwo}. Linearizing around the
standing wave solution yields a matrix potential with $V_1 =
(\beta+1)\phi^{2\beta}$ and $V_2 = \beta \phi^{2\beta}$. Hence
$V_1>0$ and $V_1>|V_2|$, which is the same as Assumption~A1).
Moreover, $L_- = -\Delta+\alpha^2 - \phi^{2\beta}$ satisfies
$L_-\phi=0$ and $L_-\ge0$ follows from $\phi>0$.

There is a large body of literature concerning the orbital (or
Lyapunov) stability (or instability) of this ground state standing
wave, see for example Shatah~\cite{Sha}, Shatah,
Strauss~\cite{ShaStr}, Weinstein~\cite{Wei1}, \cite{Wei2}, Cazenave,
Lions~\cite{CazLio}, Grillakis, Shatah, Strauss~\cite{Grill},
\cite{GSS1}, \cite{GSS2}, and Comech, Pelinovsky~\cite{ComPel}.
Reviews of much of this work are in Strauss~\cite{Str2}, and Sulem,
Sulem~\cite{SulSul}.

The question of when the stronger property of asymptotic stability
holds has received a lot of attention over the past decade. Starting
with Soffer and Weinstein~\cite{SofWei1}, \cite{SofWei2},  who
studied the modulation equations governing the evolution of {\em
small} solitons\footnote{Such solitons only arise in an NLS equation
with a linear potential. They are are generated by bifurcation off a
bound state of the linear Schr\"odinger operator.}, there has been
much work also on the case of large solitons, see Buslaev,
Perelman~\cite{BP1}, \cite{BP2}, Cuccagna~\cite{Cuc},
Perelman~\cite{Per1}, \cite{Per2}, Rodnianski, Soffer,
Schlag~\cite{RSS1}, \cite{RSS2}. It is for this purpose, rather than
for the aforementioned orbital stability, that the dispersive
estimates of the present paper are of relevance. Let us note that
for the case of small solitons the potentials $V_1, V_2$ will be
small and therefore the matrix operator above becomes easier to
treat (this is because of dimension three and analogous to the case
of scalar Schr\"odinger operators with small potentials, see e.g.,
Rodnianski, Schlag~\cite{RodSch}). Only for large $V_1, V_2$ can
significant (spectral) difficulties arise on the {\em linear} level.

It is known that Assumption~A2) implies that the spectrum
$\spec(\calH)$ satisfies $\spec(\calH) \subset \R \cup i \R$ and
that all points of the discrete spectrum other than zero are
eigenvalues whose geometric and algebraic multiplicities coincide.
For this see Grillakis~\cite{Grill}, \cite{BP1} or~\cite{RSS1}, as
well as Section~\ref{sec:gen} below.

Unfortunately it is unknown at this point how to guarantee
Assumption~A4), although it is believed to hold for systems that
arise from a ground state soliton as explained above (in 1-d this is
known, see Perelman~\cite{Per1}, due to the explicit form of the
ground state in that case). It would be desirable to have an
"abstract" approach to this question. But sofar this is unknown, and
it is an important open problem to settle this issue (even for
radial potentials). Note that there can be imbedded eigenvalues for
$V_2=0$ and $V_1$ large and positive. But in that case
Assumption~A2) does not hold. However, Assumptions A2) and A3) alone
do not imply A4) by an example\footnote{His example is in one
dimension. However, since  conditions A1)-A4) are "abstract" and
dimension less, this is relevant to our discussion.} of
Denissov~\cite{Den}. Let us remark that because of these examples
where imbedded eigenvalues can exist for our systems even though the
potentials are smooth and decay rapidly, it seems certain that the
methods known for the scalar case (say, commutator methods in the
spirit of Mourre theory) alone will not suffice. Some extra
information needs to be used (like A2 plus additional restrictions)
to insure the absence of imbedded eigenvalues.

For the case of scalar Schr\"odinger operators it is widely known
that imbedded eigenvalues are unstable. In fact, under generic
perturbations they turn into resonances in the complex plane (Fermi
golden rule). Hence, one may hope that A4) holds generically in a
suitable sense. However, in the matrix case the situation is more
complicated and imbedded eigenvalues can turn into complex
eigenvalues under small perturbations, see Cuccagna, Pelinovsky, and
Vougalter~\cite{CPV}, \cite{CP}, as well as Gang, Sigal,
Vougalter~\cite{GSV}. More precisely, whether or not this happens
depends on the sign of $\la \sigma_3 \Hil f,f\ra$ where $f$ belongs
to the real subspace associated with an imbedded eigenvalue. This is
analogous to Krein's theorem and the Krein signature in classical
mechanics, see MacKay~\cite{MR}, or Avez, Arnold~\cite{AA}.

Unlike the self-adjoint case, for our matrix operators $\Hil$ the
boundedness of $\|e^{it\Hil}\|_{2\to2}$ as $|t|\to\infty$ is
generally false. Indeed, this is the case in the presence of any
complex spectrum. Moreover, even if there is no complex spectrum,
then this operator norm can grow polynomially in $t$ due a nonzero
nilpotent part of the root-space of $\Hil$ at zero. Thus, we are
lead to consider the boundedness of $\|e^{it\Hil}P_s\|_{2\to2}$,
where $I-P_s$ is the Riesz projection corresponding to the discrete
spectrum. This has been studied before in the case where the
thresholds $\pm\mu$ are neither eigenvalues nor resonances, see
\cite{Cuc, CPV, RSS1}. In fact, the first results on such $L^2$ (or
$H^1$)-boundedness are due to Weinstein~\cite{Wei1}, \cite{Wei2} who
used variational methods. Such an approach is intimately tied up
with the underlying nonlinear problem because it uses the properties
of the ground state. For this reason, Weinstein needs to assume that
he is in the stable ($L^2$-subcritical) case. However, the recent
work~\cite{Sch1} requires such bounds also in the super-critical
case.

Our first result establishes such an $L^2$ bound in the full
generality of Assumptions A1)-A4). In particular, it shows that
neither threshold resonances nor threshold eigenvalues affect the
$L^2$-boundedness.

\begin{theorem}\label{L2}
Assume that $V$ satisfies Assumptions A1)--A4) with $\beta>5$.
Then
\[ \sup_{t\in\R}\|e^{it\Hil}P_s\|_{2\to2} \le C\]
with a constant that depends on $V$.
\end{theorem}

In this context we would like to mention the work of Gesztesy,
Jones, Latushkin, and Stanislavova~\cite{GJLS}. They prove, for
linearized NLS, that $\sigma(e^{i\Hil}P_s)= \{z\::\:|z|=1\}$.

In order to formulate our main dispersive estimate, we need to
introduce the analogue of the projection onto the continuous
spectrum from the self-adjoint case. This is done as follows.
First, let $P_d$ be the Riesz projection corresponding to the
discrete spectrum of $\Hil$. Second, let $P_\mu$ be the projection
with range equal to $\ker(\Hil-\mu)$ and kernel equal to
$(\ker(\Hil^*-\mu))^\perp$. Moreover, $P_{\mu}=0$ if $\mu$ is not
an eigenvalue of $\Hil$. Similarly with $P_{-\mu}$. We show below,
see Lemma~\ref{lem:P0}, that $P_{\pm\mu}$ are well-defined, and
that $P_d, P_\mu, P_{-\mu}$ commute. In fact, $P_d P_\mu=P_d
P_{-\mu}=P_\mu P_{-\mu}=0$.  Now, define
\[ P_c = (I-P_d)(I-P_{\mu})(I-P_{-\mu})=I-P_d-P_{-\mu}-P_\mu.\]
Clearly, $P_c$ is the analogue of the continuous spectral
projection in the self-adjoint case. It eliminates all the
eigenfunctions, including those at the thresholds (recall that we
are assuming absence of imbedded eigenvalues).

\begin{theorem}
\label{T:scalar2} Assume that $V$ satisfies Assumptions A1)--A4)
with $\beta>10$. Then there exists a time-dependent operator $F_t$
such that
$$
\sup_t\left\|F_t\right\|_{L^1\rightarrow L^\infty}<\infty, \qquad
\left\|e^{it\Hil} P_{c}-t^{-1/2}F_t
\right\|_{1\rightarrow\infty}\leq C  t^{-3/2}.
$$
If  both $\mu$ and $-\mu$ are not eigenvalues, then $F_t$ is of
rank at most two. Moreover, if $\pm\mu$ are neither eigenvalues
nor resonances, then $F_t\equiv0$.
\end{theorem}

In all cases, the operators $F_t$ can be given explicitly, and
they can  be extracted from our proofs with more work. We carry
this out explicitly for the case when $\pm \mu$ are not
eigenvalues, see formula~\eqref{ft_s2=0} below.  For scalar
Schr\"odinger operators, such explicit representations of the
kernels of $F_t$ (in terms of resonance functions and projections
onto the eigenspaces) were derived by Yajima~\cite{Y2}. His
formulas show that $F_t$ has finite rank in all cases, and the
same should be true in Theorem~\ref{T:scalar2}. It is important to
note that the $t^{-\frac32}$ bound is destroyed by an eigenvalue
at zero, even if zero is not a resonance and even after projecting
the zero eigenfunction away (this was discovered by Jensen,
Kato~\cite{JenKat} for scalar operators).

Finally, we remark that it was not our intention to obtain the
minimal value of $\beta$ in Assumption~A3). Our results can surely
be improved in that regard. Needless to say, the problem of lowering
the requirement on $\beta$ is only one of many remaining issues.
More relevant to nonlinear questions seems to be how to prove A4),
and/or how to deal with imbedded eigenvalues when they do occur (in
regards to our theorems). In a similar vein, it would of course be
interesting to develop this linear theory when A2) does not hold.
This is the case, for example, when linearizing around excited
states, see~\cite{BL2}.

\section{The matrix case: Generalities}
\label{sec:gen}

In this section we shall develop some standard and well-known
properties of the spectra and resolvents of $\Hil$ under
Assumptions A2)-A4). It should be mentioned that Assumption A1)
seems to be needed only in order to apply the symmetric resolvent
identity, see Section~\ref{sec:thresh} below. However, in this
section we work with the usual resolvent identity and therefore do
not need~A1).\footnote{It seems that one can work with the usual
resolvent identity throughout this paper, which would then allow
us to dispense with A1) altogether. However, A1) holds in
important applications and we find it convenient to use the
symmetric resolvent identity.}

\begin{lemma}
\label{lem:spec} Let $\beta>0$ be arbitrary in~\eqref{decay}. Then
the  essential spectrum of $\Hil$ equals
$(-\infty,-\mu]\cup[\mu,\infty)$. Moreover,
$\spec(\Hil)=-\spec(\Hil)=\overline{\spec(\Hil)}=\spec(\Hil^*)$
and $\spec(\Hil)\subset \R\cup i\R$. The discrete spectrum of
$\Hil$ consists of eigenvalues $\{z_j\}_{j=1}^N$, $0\le N\le
\infty$, of finite multiplicity. For each $z_j\ne0$ the algebraic
and geometric multiplicities coincide and $\Ran(\Hil-z_j)$ is
closed. The zero eigenvalue has finite algebraic multiplicity,
i.e., the generalized eigenspace $\bigcup_{k=1}^\infty
\ker(\Hil^k)$ has finite dimension. In fact, there is a finite
$m\ge1$ such that $\ker(\Hil^k)=\ker(\Hil^{k+1})$ for all $k\ge
m$.
\end{lemma}
\begin{proof}
The statement about the essential spectrum follows from Weyl's
criterium. To see this, note that conjugation of $\Hil$ by the
matrix $\bm 1 & i \\ 1 & -i\endm$ leads to the matrix operator
\[ \bm 0 & i L_{-} \\ -iL_{+} & 0 \endm \]
where $L_-$ is as above and with $L_+ = - \Laplace + \mu - V_1 -
V_2$. We will again denote this matrix by $\Hil$. Let
$H_0=-\Laplace+\mu$ and set $W_1=-V_1+V_2,\;W_2=-V_1-V_2$,
\begin{align}
 \label{eq:Hildef}
 \Hil_0 &= \bm 0 &i H_0 \\ -i H_0 & 0 \endm,\quad W= \bm 0 & iW_1 \\ -iW_2 & 0
\endm,\\
  \Hil &= \Hil_0+W=i\bm 0 & H_0+W_1 \\ -H_0-W_2 & 0 \endm. \nn
  \end{align}
By means of the matrix $J=\bm 0 & i \\ -i & 0 \endm$ one can also
write
\[ \Hil_0 = \bm H_0 & 0 \\ 0 & H_0 \endm J,\quad \Hil = \bm H_0+W_1 & 0 \\ 0 &
H_0+W_2 \endm J.
\]
Clearly, $\Hil$ is a closed operator on the domain
$\Dom(\Hil)=W^{2,2}\times W^{2,2}$. Since $\Hil_0^*=\Hil_0$ it
follows that $\spec(\Hil_0)\subset \R$. One checks that for $\Re
z\ne0$ \begin{align} (\Hil_0-z)^{-1} &= -(\Hil_0+z)\bm
(H_0^2-z^2)^{-1} & 0
\\ 0 & (H_0^2-z^2)^{-1}
\endm  \nn \\
&= -\bm (H_0^2-z^2)^{-1} & 0 \\ 0 & (H_0^2-z^2)^{-1} \endm
(\Hil_0+z)
\label{eq:Bspec} \\
(\Hil-z)^{-1} &= (\Hil_0-z)^{-1}\!\! -\! (\Hil_0-z)^{-1} U_1
\Bigl[1+U_2 J (\Hil_0-z)^{-1} U_1\Bigr]^{-1}\! U_2 J
(\Hil_0-z)^{-1} \label{eq:grill} \end{align} where
\eqref{eq:grill} also requires the expression in brackets to be
invertible, and with
\[
 U_1 = \bm |W_1|^{\half} & 0 \\ 0 & |W_2|^{\half} \endm,\quad U_2 =
\bm |W_1|^{\half}\sign(W_1) & 0 \\ 0 & |W_2|^{\half}\sign(W_2)
\endm.
\]
It follows from \eqref{eq:Bspec} that $\spec(\Hil_0) =
(-\infty,-\mu]\cup [\mu, \infty)\subset \R$. Since $V_1(x)\to0$
and $V_2(x)\to0$ as $x\to\infty$, it follows from Weyl's theorem,
see Theorem~XIII.14 in~\cite{RS4}, and the
representation~\eqref{eq:grill} for the resolvent of~$\Hil$, that
$\spec_{ess}(\Hil)=\spec_{ess}(\Hil_0)=(-\infty,-\mu]\cup [\mu,
\infty)\subset \R$. Moreover, \eqref{eq:grill} implies via the
analytic Fredholm alternative that $(\Hil-z)^{-1}$ is a
meromorphic function in~$\C\setminus (-\infty,-\mu]\cup [\mu,
\infty)$.  Furthermore, the poles are eigenvalues\footnote{Note
that since $\Hil$ is not self-adjoint, it can happen that
\[ \ker(\Hil-z)^2 \ne \ker (\Hil-z)\]
for some $z\in\C$. In other words, $\Hil$ can possess {\em
generalized eigenspaces}. In the NLS applications this does happen
at $z=0$ due to symmetries like modulation.} of~$\Hil$ of finite
multiplicity and $\Ran(\Hil-z_j)$ is closed at each pole $z_j$.

The symmetries of the spectrum are consequences of the commutation
properties of $\Hil$ with the Pauli matrices
\[ \sigma_1=\bm 0 & 1\\ 1& 0\endm, \quad \sigma_2=\bm 0 & i\\ -i &  0\endm,
\quad \sigma_3=\bm 1 & 0\\ 0 & -1\endm. \] Now let us check that
the spectrum lies in the union of the real and imaginary axes.
Thus, suppose that
\[ \bm 0 & i L_{-} \\ -iL_{+} & 0 \endm \binom{f_1}{f_2} = E
\binom{f_1}{f_2} \] with $E\ne0$ and $\binom{f_1}{f_2}\in
L^2\setminus\{0\}$. Then $f_1\ne0$ and $f_2\ne0$ and $f_1\perp
\ker\{L_-\}$. Hence,  $g= L_-^{-\half}f_1$ satisfies
\[ L_-^{\half} L_+ L_-^{\half} g = E^2 g \]
and thus $E^2 \in \R$, as desired. Here we used that
$\sqrt{L_{-}}L_{+}\sqrt{L_{-}}$ with domain $W^{4,2}(\R^3)$ is a
selfadjoint operator. For a proof of this see Lemma~11.10
in~\cite{RSS2}. That same lemma also contains a proof of the fact
that for any eigenvalue other than zero the algebraic and
geometric multiplicities coincide.

Let $P_0$ be the Riesz projection at zero. Then, on the one hand
one checks that
\[ \Ran P_0 \supset \ker(\Hil^m) \text{\ \ for all\ \ }m\ge1.\]
On the other hand, if $\|(\Hil-z)^{-1}\|\les |z|^{-\nu}$, then
\[ \Hil^\nu P_0 =0.\]
Thus $\Ran P_0\subset \ker (\Hil^\nu)$. See~\cite{HS} Chapter~6
for these general statements about Riesz projections.
\end{proof}

It will follow from the next section that $N<\infty$ in
Lemma~\ref{lem:spec} provided $\beta>5$ (which can probably be
relaxed). Indeed, in that section we will derive expansions of the
resolvent $(\Hil-z)^{-1}$ about the thresholds $\pm\mu$ which will
preclude the eigenvalues from accumulating at these points. Thus
there can only be finitely many eigenvalues, i.e., $N<\infty$.

Next, we need to develop a limiting absorption principle for the
resolvents $(\Hil-z)^{-1}$ when $|z|>\mu$. As observed in \cite{CPV}
and \cite{Sch1}, this can be done along the lines of the classical
Agmon argument~\cite{agmon}. We now present some of these arguments.

We begin by recalling some weighted $L^2$ estimates for the free
resolvent $(\Hil_0-z)^{-1}$ which go by the name "limiting
absorption principle". The weighted $L^2$-spaces here are the
usual ones $L^{2,\sigma}=\la x\ra^{-\sigma} L^2$. It will be
convenient to introduce the space
\[ X_\sigma := L^{2,\sigma}(\R^3)\times L^{2,\sigma}(\R^3). \]
Clearly, $X_\sigma^* = X_{-\sigma}$. The statement is that
\beeq
\label{eq:H0agm}
 \sup_{|\lambda|\ge\lambda_0,\,0<\eps} |\lambda|^{\half}\,
 \|(\Hil_0-(\lambda\pm i\eps))^{-1}\|_{X_\sigma\to X_\sigma^*} < \infty
\eneq provided $\lambda_0>\mu$ and $\sigma>\half$ and was proved
in this form by Agmon~\cite{agmon}. By the explicit expression for
the kernel of the free resolvent in $\R^3$ one obtains the
existence of the limit
\[ \lim_{\eps\to0+} \la (\Hil_0-(\lambda\pm i\eps))^{-1} \phi,\psi \ra \]
for any $\lambda\in\R$ and any pair of Schwartz functions
$\phi,\psi$. Hence $(\Hil_0-(\lambda\pm i0))^{-1}$  satisfies the
same bound as in~\eqref{eq:H0agm} provided
$|\lambda|\ge\lambda_0>\mu$. There is a corresponding  bound which
is valid for all energies. It takes the form \beeq \label{eq:HS1}
 \sup_{z\in\C}\|(\Hil_0-z)^{-1}\|_{X_\sigma\to X_\sigma^*} < \infty
\eneq provided $\sigma>1$. It is much more elementary  to obtain
than~\eqref{eq:H0agm} since it only uses that the convolution with
$|x|^{-1}$ is bounded from $L^{2,\sigma}(\R^3)\to
L^{2,-\sigma}(\R^3)$ provided $\sigma>1$. In fact, it is
Hilbert-Schmidt in these norms. We now state a lemma about absence
of imbedded resonances.

\begin{lemma}
\label{lem:invert} Let $\beta>1$. Then for any $\lambda\in\R$,
$|\lambda|>\mu$ the operator $(\Hil_0-(\lambda\pm i0))^{-1}V$ is a
compact operator on $X_{-\half-}\to X_{-\half-}$ and
\[ I + (\Hil_0-(\lambda\pm i0))^{-1}V \]
is invertible on these spaces.
\end{lemma}
\begin{proof}
The compactness is standard and we refer the reader
to~\cite{agmon} or~\cite{RS4}. Let $\lambda>\mu$. By the Fredholm
alternative, the invertibility statement requires excluding
solutions $(\psi_1,\psi_2)\in X_{-\half-}$ of the system
\begin{align*}
0 &= \psi_1 -  R_0(\lambda-\mu+i0)(V_1\psi_1 + V_2\psi_2) \nn \\
0 &= \psi_2 -  R_0(-\lambda-\mu)(V_2\psi_1+V_1\psi_2), \nn
\end{align*}
where $R_0(z)$ is the free, scalar resolvent $(-\Laplace-z)^{-1}$.
Notice that these equations imply that $\psi_2\in L^2$ and that
\begin{align*}
0 &= \la \psi_1, V_1\psi_1\ra + \la \psi_1, V_2\psi_2\ra - \la
R_0(\lambda-\mu+i0)(V_1\psi_1+V_2\psi_2),
V_1\psi_1+V_2\psi_2 \ra \nn \\
0 &= \la \psi_2, V_2\psi_1 \ra - \la R_0(-\lambda-\mu)(V_2\psi_1+V_1\psi_2), V_2\psi_1 \ra \nn \\
0 &= \la \psi_2, V_1 \psi_2 \ra - \la R_0(-\lambda-\mu) V_2\psi_1,
V_1\psi_2\ra - \la R_0(-\lambda-\mu)V_1\psi_2, V_1\psi_2\ra. \nn
\end{align*}
Since $V_1,V_2$ are real-valued, inspection of these equations
reveals that
\[ \Im \la R_0(\lambda-\mu+i0)(V_1\psi_1+V_2\psi_2), V_1\psi_1+V_2\psi_2\ra =0.\]
So Agmon's well-known bootstrap lemma (see Theorem~3.2
in~\cite{agmon}) can be used to conclude that $\psi_1\in
L^2(\R^3)$. But then we have an imbedded eigenvalue at $\lambda$,
which contradicts Assumption~A4). So one can invert
\[ I+ (\Hil_0-(\lambda\pm i0))^{-1}V\]
on $X_{-\half-}$ and we are done.
\end{proof}

As usual, one converts the information of the previous lemma into
a bound for the perturbed resolvent by means of the resolvent
identity.

\begin{proposition}
\label{prop:lim_ap} Let $\beta>1$ and fix an arbitrary
$\lambda_0>\mu$. Then \beeq \label{eq:limap}
\sup_{|\lambda|\ge\lambda_0,\,0<\eps}
|\lambda|^{\half}\,\|(\Hil-(\lambda\pm i\eps))^{-1}\| < \infty
\eneq where the norm is the one from  $X_{\half+}\to X_{-\half-}$.
\end{proposition}
\begin{proof}
 Let $z=\lambda+i\eps$, $\lambda\ge\lambda_0$,
$\eps\ne0$. By the resolvent identity and the fact that the
spectrum of $\Hil$ belongs to $\R\cup i\R$,
\beeq \label{eq:res}
 (\Hil-z)^{-1} = (I+(\Hil_0-z)^{-1}V)^{-1}(\Hil_0-z)^{-1}
\eneq
as operators on $L^2(\R^3)$. Because of the
$|\lambda|^{-\half}$-decay in~\eqref{eq:H0agm}, there exists a
positive radius $r_V$ such that
\[ \|(\Hil_0-z)^{-1}V\|<\half \]
for all $|z|>r_V$ in the operator norm of $X_{-\half-}\to
X_{-\half-}$. In conjunction with~\eqref{eq:res} this implies that
\[ \|(\Hil-z)^{-1}\| \les |z|^{-\half} \]
for all $|z|>r_V$ in the operator norm of $X_{\half+}\to
X_{-\half-}$. Now suppose~\eqref{eq:limap} fails. It then follows
from~\eqref{eq:res} and~\eqref{eq:HS1} that there exist a sequence
$z_n$ with $\Re(z_n)\ge\lambda_0$ and functions $f_n \in
X_{-\half-}$ with $\|f_n\|_{X_{-\half-}}=1$ and such that \beeq
\label{eq:small_norm}
 \|[I+(\Hil_0-z_n)^{-1}V] f_n \|_{X_{-\half-}} \to 0
\eneq as $n\to\infty$. Necessarily, the $z_n$ accumulate at some
point $\lambda\in [\lambda_0,r_V]$. Without loss of generality,
$z_n\to \lambda$ and $\Im(z_n)>0$ for all $n\ge1$. Next, we claim
that~\eqref{eq:small_norm} also holds in the following form:
\beeq
\label{eq:small_norm2}
 \|[I+(\Hil_0-(\lambda+i0))^{-1}V] f_n \|_{X_{-\half-}} \to 0
\eneq as $n\to\infty$. If so, then it would clearly contradict
Lemma~\ref{lem:invert}.  To prove~\eqref{eq:small_norm2}, let
\[ S:= I+(\Hil_0-(\lambda+i0))^{-1}V \]
for simplicity. Then
\begin{align}
I+(\Hil_0-z_n)^{-1}V &= S + ((\Hil_0-z_n)^{-1}-(\Hil_0-(\lambda+i0))^{-1}) V \nn \\
&= \big[I + ((\Hil_0-z_n)^{-1}-(\Hil_0-(\lambda+i0))^{-1})
VS^{-1}\big]S. \label{eq:SHnew}
\end{align}
Our claim now follows from the fact that the expression in
brackets is an invertible operator for large $n$ on $X_{-\half-}$.
This in turn relies on bounds of the form: Given $\eps>0$, there
exists $\delta>0$ so that  for $\Re z>0$, and all $z'$ close to $
z$, \beeq \label{eq:ag_hold} \big\|
(-\Delta-z)^{-1}-(-\Delta-z')^{-1} \big\|_{L^{2,\half+\eps}\to
L^{2,-\half-\eps}} \le C_{\delta,\eps}\,|z-z'|^{\delta} \eneq
see~\cite{agmon}.\footnote{Of course $\delta\to0$ as $\eps\to0$.
Moreover, if $\delta=1$, then  one needs $\eps>1$.}
\end{proof}

As in the case of the free Hamiltonian $\Hil_0$, it is now
possible to define the boundary values of the resolvent
$(\Hil-z)^{-1}$. More precisely, the following corollary holds.

\begin{cor}
\label{cor:Hbdry} Let $\beta>1$. Define
 \beeq \label{eq:Hbdrydef}
(\Hil-(\lambda\pm i0))^{-1}:=(I+(\Hil_0-(\lambda\pm
i0))^{-1}V)^{-1}(\Hil_0-(\lambda\pm i0))^{-1} \eneq for all
$|\lambda|>\mu$.  Then as $\eps\to0+$,
\[ \|(\Hil-(\lambda\pm i\eps))^{-1} - (\Hil-(\lambda\pm i0))^{-1}\| \to0\]
in the norm of $X_{\half+}\to X_{-\half-}$ and one can extend
\eqref{eq:limap} to $\eps\ge0$.
\end{cor}
\begin{proof} Definition \eqref{eq:Hbdrydef}
is legitimate by Lemma~\ref{lem:invert} and motivated
by~\eqref{eq:res}. Thus, the resolvent $(\Hil-(\lambda\pm
i\eps))^{-1}$ is well-defined for all $\eps\ge0$ and
$|\lambda|>\lambda_0$. In view of~\eqref{eq:ag_hold},
\[  \|(\Hil_0-(\lambda\pm i\eps))^{-1} - (\Hil_0-(\lambda\pm i0))^{-1}\| \to 0\]
as $\eps\to0$ in the norm of $X_{\half+}\to X_{-\half-}$.
Moreover, by~\eqref{eq:SHnew} and again~\eqref{eq:ag_hold},
\begin{align*}
& [I+(\Hil_0-(\lambda+i\eps))^{-1}V]^{-1} - [I+(\Hil_0-(\lambda+i0))^{-1}V]^{-1} \\
& = S^{-1}\big[I + ((\Hil_0-(\lambda+i\eps))^{-1}-(\Hil_0-(\lambda+i0))^{-1}) VS^{-1}\big]^{-1} - S^{-1} \\
& = \sum_{k=1}^\infty S^{-1}
\big[-((\Hil_0-(\lambda+i\eps))^{-1}-(\Hil_0-(\lambda+i0))^{-1})
VS^{-1}\big]^{k}
\end{align*}
tends to zero in the norm of $X_{-\half-}$ as $\eps\to0+$.
\end{proof}

\section{Resolvent expansions at thresholds}
\label{sec:thresh}

In view of Assumption~A1), we write
$$
V = - \sigma_3 v v^* = -\sigma_3 v^2=:v_1v_2,
$$
where $v_1 = - \sigma_3 v,$  $v_2 = v^* = v,$
\begin{equation}
\label{eq:pauli}
\sigma_3=\left[\begin{array}{cc} 1 & 0 \\ 0 & -1
\end{array}\right].
\end{equation}
It follows from~\eqref{decay} that the entries of $v_1, v_2$ are
real-valued and decay like $\la x\ra^{-\beta/2}$.  Let $\lambda =
\mu + z^2$, where $\Ima(z) > 0$ and $|z|$ small. We have the
symmetric resolvent identity: \bea \label{mresolve}
R(\lambda):=(\calH - \lambda)^{-1} &= R_0(\lambda) - R_0(\lambda)
v_1(I+v_2 R_0 (\lambda) v_1)^{-1} v_2 R_0(\lambda).
\eea

Recall that (see previous section)  the
essential spectrum of $\calH$ is $(-\infty, -\mu] \cup
[\mu,\infty)$. As in the scalar case \cite{ErdSch}, we obtain resolvent expansions at the threshold
$\lambda=\mu$ in the case of a resonance and/or eigenvalue.
Recall that $R_0(\lambda)$ has the kernel
$$
R_0 (\lambda)(x,y)=\f{1}{4\pi|x-y|}\left[
\begin{array}{cc}
 e^{iz|x-y|} & 0\\
0&- e^{-\sqrt{2\mu+z^2}|x-y|}
\end{array}
\right].
$$
We have a similar representation of $R_0 (\lambda)$ for $\lambda$ around $-\mu$.
Let
\begin{align*}
A(z) &= I + v_2 R_0(\lambda) v_1\\
&=: A_0+zA_1(z),
\end{align*}
where
\begin{align*}
A_0 &= I + v_2 R_0 (\mu) v_1,\\
A_1(z)&=: \f{1}{z}v_2\left(R_0(\lambda)-R_0(\mu)\right)  v_1, \\
R_0(\mu)(x,y)&=\f{1}{4\pi|x-y|}\left[
\begin{array}{cc}
 1 & 0\\
0&- e^{-\sqrt{2\mu}|x-y|}
\end{array}
\right].
\end{align*}
If $\beta>3$, then $v_2 R_0(\lambda) v_1$ is a Hilbert-Schmidt
operator. Hence, $\ker A_0$ is finite-dimensional.

\begin{lemma}\cite{jensennenciu}\label{L:jen-nen}
Let $F\subset \C\setminus\{0\}$ have zero as an accumulation point. Let $A(z)$, $z\in F$,
be a family of bounded operators of the form
$$A(z)=A_0+z A_1(z)$$
with $A_1(z)$ uniformly bounded as $z\rightarrow 0$. Suppose that $0$ is an isolated point of the spectrum of
$A_0$, and let $S$ be the corresponding Riesz projection. Assume that $\rank(S)<\infty$.
Then for sufficiently small $z\in F$ the operators
$$
B(z):=\frac{1}{z}(S-S(A(z)+S)^{-1} S)
$$
are well-defined and bounded on $\calH$. Moreover, if $A_0=A_0^*$, then they are
 uniformly bounded as $z\rightarrow 0$. The operator $A(z)$ has a bounded inverse in $\calH$ if and only if
$B(z)$ has a bounded inverse in $S\calH$, and in this case
\bea\label{az-1}
A(z)^{-1}=(A(z)+S)^{-1}+\frac{1}{z}(A(z)+S)^{-1}SB(z)^{-1}S(A(z)+S)^{-1}.
\eea
\end{lemma}
\noindent See \cite{ErdSch} for the proof.

We use Lemma~\ref{L:jen-nen} to obtain an expansion of
$A(z)^{-1}$. Assume $A_0$ is not invertible. Let $S_1$ be the
Riesz projection corresponding to $0$. As in the scalar case,
$A_0$ is self adjoint and it is a compact perturbation of the
identity. Therefore, $S_1 = P_{\ker A_0}$, $A_0+S_1$ is invertible
and \bea \label{s1commute}
S_1=(A_0+S_1)^{-1}S_1=S_1(A_0+S_1)^{-1}. \eea Also note that, if
$V$ satisfies \pref{decay} for some $\beta>3$, then \bea \nonumber
\sup_{\begin{subarray}{c}
\text{$z$ small}\\
\text{Im}(z) \geq 0
\end{subarray}}
\| A_1(z) \|_{HS}<\infty
\eea
Thus, $A(z) + S_1$ is invertible for small $z$.
By Lemma~\ref{L:jen-nen} we have
\begin{align*}
A(z)^{-1} &= \left(A(z) +S_1\right)^{-1} + \f{1}{z}\left(A(z) + S_1\right)^{-1} S_1
m(z)^{-1} S_1 \left(A(z) + S_1\right)^{-1},
\end{align*}
where
\begin{align*}
m(z) &= \f{1}{z} \left(S_1 - S_1 (A(z)+S_1)^{-1}S_1\right)\\
&=  \f{-1}{z} S_1 \left[
\sum_{k=1}^\infty (-1)^k z^k \left(A_1(z) (A_0 + S_1)^{-1}\right)^k \right] S_1\\
&= S_1 A_1(z) S_1 +   \sum_{k=1}^\infty
(-1)^k z^k  S_1\left(A_1(z) (A_0 + S_1)^{-1}\right)^{k+1}S_1\\
&=: S_1 A_1 (0) S_1 + zm_1(z).
\end{align*}
We used (\ref{s1commute}) in the second equality.
Let $f=\binom{f_1}{f_2}$ and define
$$
P_1f:=\int_{\R^3}\left[\begin{array}{cc}
1&0\\
0&0
\end{array}
\right]f(x)\, dx
\text{ and }
P_2f:=\int_{\R^3}\left[\begin{array}{cc}
0&0\\
0&1
\end{array}
\right]f(x)\, dx.
$$
Note that
$$
A_1(0) = \f{i}{4\pi} v_2 P_1 v_1.
$$
Therefore,
$$
m(0) = \f{-i}{4\pi} S_1 v P_1 v S_1.
$$
As in the scalar case, if $m(0)$ is invertible, then we invert $m(z)$ using Neumann series. Otherwise, let
$S_2 = P_{\ker m(0)}: S_1 L^2 \to S_1 L^2$.
Obviously, $m(z) + S_2$ is invertible for small $z$. Using Lemma~\ref{L:jen-nen}, we have
$$
m(z)^{-1} =  \left(m(z) + S_2 \right)^{-1} + \f{1}{z}\left(m(z) + S_2 \right)^{-1}
S_2  b(z)^{-1} S_2 \left(m(z) + S_2\right)^{-1},
$$
where
\begin{align} \label{eq:bdef}
b(z) &= \f{1}{z}\left(S_2 - S_2 \left(m(z) + S_2\right)^{-1} S_2\right)\\
&= \f{-1}{z}  S_2 \left[\sum_{k=1}^\infty (-1)^k z^k
\left(m_1(z)(m(0)+S_2)^{-1}\right)^k\right]S_2  \nn \\
 &=:b(0)+ zb_1(z), \nn
\end{align}
where
$$
b(0)=S_2 m_1 (0) S_2.
$$
Note that \beq \label{eq:m1def}
m_1(z)=S_1\f{A_1(z)-A_1(0)}{z}S_1+\sum_{k=1}^\infty (-1)^k z^{k-1}
S_1\left(A_1(z) (A_0 + S_1)^{-1}\right)^{k+1}S_1. \eneq Therefore
\begin{align} \label{eq:b0def}
b(0)&=S_2 m_1 (0) S_2\\
&=\f{1}{8\pi}S_2v_2\left[\begin{array}{cc}
- |x-y| & 0\\
0& \f{e^{\sqrt{2\mu}|x-y|}}{\sqrt{2\mu}}
\end{array}
\right]
v_1 S_2 \nn \\
&= \f{1}{8\pi} S_2 v \left[\begin{array}{cc}
|x-y|&0\\
0&\f{1}{\sqrt{2\mu}}e^{-\sqrt{2\mu}|x-y|}\\
\end{array}\right]
vS_2. \nn
\end{align}

Below, we will characterize the projections $S_1$, $S_2$ and prove that $b(0)$ is always invertible in $S_2L^2$.
\begin{lemma}\label{L:ms1s2}
Assume $\beta>3$. Then
\\
i) $f \in S_1L^2\backslash\{0\}$ if and only if $f=v_2g$ for some $g \in L^{2,-\f{1}{2}-}\backslash\{0\}$ such that
\bea \label{mHg=0}
(\calH_0-\mu)g + Vg =0 \ \text{in}\ \calS^\prime.
\eea
ii) Assume $f\in S_1L^2\backslash\{0\}$, then the following are equivalent\\
a) $f \in S_2L^2\backslash\{0\}$,\\
b) $P_1vf=0$,\\
c) $f=v_2g$ for some $g \in L^{2}\backslash\{0\}$ satisfying \pref{mHg=0}.
\end{lemma}
\begin{proof}
If $f \in S_1L^2\backslash\{0\}$, then
\[ A_0f = f + v_2 R_0(\mu) v_1 f =0 \]
by definition. Hence, $f=v_2 g $ where
\[ g = - R_0(\mu) v_1 f \in L^{2,-\f{1}{2}-} (\R^3)\]
by the mapping properties of $(-\Delta)^{-1}$. Moreover,
\[ g + R_0(\mu) V g =0. \]
 By Lemma~2.4 in~\cite{JenKat} this is equivalent
with~\eqref{mHg=0}. Conversely, if~\eqref{mHg=0} holds, then we
set $f=v_2g$ which belongs to $L^2$ and satisfies
\[ f + v_2 R_0(\mu) v_1 f =0. \]
Thus, $S_1 f =f$, and the first part is proven.

For the second part, suppose that $S_1f=f$. If in addition
$S_2f=f$, then $m(0)f=0$ which is the same as $S_1 vP_1 vf=0$. But
then also
\[ \la S_1vP_1 vf,f\ra = \la P_1 vf,vf \ra = \binom{(P_1 vf)_1^2}{0}
\]
where we have written $P_1vf = \binom{(P_1 vf)_1}{0}$. Hence
$P_1vf=0$. This
implies that
\[ g = - R_0(\mu) v_1 f \in L^{2} (\R^3)\]
This is a standard property, see for example Lemma~6
in~\cite{ErdSch}. In view of the first part of this proof
$f=v_2g$.

These implications can be reversed: Indeed, if
\[ g = - R_0(\mu) v_1 f \in L^{2} (\R^3)\]
then it follows easily that $P_1 v_1f=0$ which is the same as $P_1
vf=0$ (see for example Lemma~6 in~\cite{ErdSch}). But then also
$m(0)f=0$, and the lemma follows.
\end{proof}

Next, we show that the Jensen-Nenciu expansion stops after (at
most) two steps.

\begin{lemma}  \label{L:mb0}
Assume $\beta>5$. Then, as an operator in $S_2L^2$, the kernel of
$b(0)$ is trivial.
\end{lemma}
\begin{proof}
Assume $f \in S_2L^2$ is in $\ker b(0)$. Since $b(0)$ has a
real-valued kernel, we can assume that $f$ is real-valued. Let
$f=\binom{f_1}{f_2}$ and $h= vf =\binom{h_1}{h_2}$. By
Lemma~\ref{L:ms1s2} ii), we have $\int h_1 =0$, $h\equiv -\sigma_3
Vg$ for some real-valued $g=\binom{g_1}{g_2}\in L^{2,-\half-}$
satisfying (\ref{mHg=0}) and
\bea \label{hg} h_1 = V_1 g_1 + V_2
g_2 \;, \ h_2 = V_2 g_1 + V_1 g_2. \eea
Moreover, since $f\in \ker
b(0)$ (again by Lemma~\ref{L:ms1s2}), we have \bea \label{hh}
\left\langle h \;, \left[\begin{array}{cc}
|x-y| & 0\\
0 & \f{1}{\sqrt{2\mu}} e^{-{\sqrt{2\mu}} |x-y|}
\end{array}
\right] h\right\rangle = 0. \eea
Now use the following fact
from~\cite{JenKat} (see also the proof of Lemma~7
in~\cite{ErdSch}): if $\int u = \int v = 0$, and $u,v \in
L^{2,s}$, $s>5/2$, then
\begin{align*}
\langle |x-y| u,v\rangle &= -\f{1}{2\pi}
\Big\langle \f{1}{|x|} \ast u \;, \ \f{1}{|x|} \ast v\Big\rangle\\
&= - 8\pi \langle ( - \Delta )^{-1} u\;, \ (-\Delta)^{-1} v \rangle.
\end{align*}
Thus,
\begin{align*}
\eqref{hh} & = - 8\pi \| (-\Delta)^{-1} h_1 \|_2^2 +
\langle h_2 \;, \ \f{1}{\sqrt{2\mu}}
e^{-\sqrt{2\mu}|x-y|} h_2 \rangle =0
\end{align*}
Define $\hat{f} (\xi) = \int e^{-i\xi\cdot x} f(x) \; dx $.
  Recall that (see, e.g., \cite{Ste})
\begin{align*}
\f{\widehat{e^{-\sqrt{2\mu}|x|}}}{4\pi|x|}(\xi) &= (\xi^2 + 2\mu)^{-1},\\
\f{\widehat{e^{-\sqrt{2\mu}|x|}}}{\sqrt{2\mu}}(\xi) &=
\f{8\pi}{(\xi^2 + 2\mu)^2}.
\end{align*}
Thus,
\bea \label{g1=g2}
  \| (-\Delta )^{-1} h_1 \|^2_2 =
  \|(-\Delta + 2\mu)^{-1}h_2 \|^2_2.
\eea
On the other hand, by (\ref{mHg=0}), we have
\bea \label{Hg=mg}
\left[\begin{array}{cc}
-\Delta - V_1 & -V_2\\
V_2 & \Delta - 2\mu + V_1
\end{array}
\right]
\binom{g_1}{g_2} = 0.
\eea
Using this and (\ref{hg}), we obtain
\begin{align*}
-\Delta g_1 - V_1 g_1 - V_2 g_2 = 0 &\qquad\Rightarrow &h_1 &= -\Delta g_1,\\
\Delta g_2 - 2\mu g_2 + V_1g_2 + V_2 g_1 = 0&\qquad\Rightarrow &h_2 &=(-\Delta + 2\mu) g_2.
\end{align*}
Adding the equalities on the left hand side, we obtain
\[
L_- (g_1 - g_2) = (-\Delta + \mu - V_1 + V_2 )(g_1 - g_2)=\mu (g_1+ g_2).
\]
Pairing this with $g_1-g_2$, we have (recall that $g_1,g_2$ are
real-valued)
\begin{align*}
\langle L_- (g_1 - g_2) , g_1 - g_2 \rangle &= \mu
\left(\| g_1\|^2_2 - \|g_2 \|^2_2\right)\\
&= 0 \ \text{by (\ref{g1=g2})}
\end{align*}
The positivity assumption $L_- \geq 0$ implies that
$\ker L_- = \text{span} \{\varphi\}$ (if $\ker L_- = \{0\}$, then  $\varphi=0$.
Otherwise $\varphi \not=0$). Therefore,
\[
g_1 - g_2 = k\varphi \;, \; \text{ for some }k\in \R\;.
\]
Using this in (\ref{Hg=mg}), we have
\begin{align*}
\left[\begin{array}{cc}
-\Delta - V_1 & -V_2\\
V_2 & \Delta - 2\mu + V_1
\end{array}
\right]
\left(\begin{array}{c}
g_1\\
g_1 - k \varphi
\end{array}
\right) &=0 \;\;\Rightarrow \\
(-\Delta - V_1 - V_2) g_1 + k V_2 \varphi &=0,\\
(\Delta - 2\mu + V_1 + V_2) g_1 - k ( \Delta - 2\mu + V_1) \varphi &=0.
\end{align*}
Adding the last two inequalities and using the fact that $\varphi\in \ker L_-$, we have
\[
g_1 = \f{k}{2} \varphi
\ \Rightarrow
g_2 = -\f{k}{2} \varphi.
\]
If $k\ne 0$, we use (\ref{Hg=mg}) once more to conclude that
\[
\left[\begin{array}{cc}
-\Delta -V_1 & -V_2\\
V_2 & \Delta - 2\mu + V_1
\end{array}
\right]
\left(\begin{array}{c}
\phantom{-}\varphi\\
-\varphi
\end{array}
\right)=0.
\]
This implies that
\[
(-\Delta - V_1 + V_2 ) \varphi = 0 \quad
\Rightarrow\quad\mu\varphi =0 \quad \Rightarrow\quad \varphi
\equiv 0
\]
Hence, in all cases $g_1=g_2=g$. But then
\begin{align*}
-\Delta g - V_1 g - V_2 g = 0 \\
\Delta g - 2\mu g + V_1g + V_2 g = 0
\end{align*}
which implies that $\mu g=0$ and thus also $g=0$. Retracing our
steps we conclude that $h=0$ and $f=0$. Therefore, $\ker
b(0)=\{0\}$ and we are done.
\end{proof}

Lemmas \ref{L:jen-nen}, \ref{L:ms1s2} and \ref{L:mb0} imply that
$A(z)$ is always invertible for small $z\ne0$ and
\begin{align}
 &A(z)^{-1}  = \left(A(z) + S_1\right)^{-1} \label{eq:Aexp}\\
 &+ \f{1}{z} \left(A(z) + S_1\right)^{-1} S_1
\left(m(z) + S_2\right)^{-1} \left(A(z) + S_1\right)^{-1} + \nn\\
&+\f{1}{z^2} \left(A(z) + S_1\right)^{-1} S_1 \left(m(z) +
S_2\right)^{-1} S_2 b(z)^{-1} S_2 \left(m(z) + S_2\right)^{-1} S_1
\left(A(z) + S_1\right)^{-1}.\nn
\end{align}
Note that
\bea \nonumber
  A(z)^{-1}  =  \f{1}{z^2} S_2 b(0)^{-1}S_2 + O(\f{1}{z}).
\eea
With $\lambda = \mu + z^2$,
\begin{align}
R_V(\lambda) &= R_0 (\lambda ) - R_0 (\lambda ) v_1 \left(A(z)\right)^{-1} v_2 R_0 (\lambda)\nn \\
&= -\f{1}{z^2} R_0(\lambda) \sigma_3 v  S_2 b(0)^{-1}   S_2 v  R_0 (\lambda) + \ldots \label{eq:RV}
\end{align}
The most singular term in this expansion can be identified as a
(not necessarily orthogonal) projection onto the eigenspace at the
threshold.

\begin{lemma}
\label{lem:P0} Let $\beta>5$. Then the operator $P_\mu:=-R_0(\mu)
\sigma_3 v S_2 b(0)^{-1} S_2 v  R_0(\mu)$ is a projection in
$L^2(\R^3)\times L^2(\R^3)$ with the property that
\[ \ran(P_\mu)=\ker (\calH-\mu), \qquad \ker(P_\mu)= \ker (\calH^\ast-\mu)^\perp. \]
\end{lemma}
\begin{proof}
Choose  a basis $\{\varphi_j\}_{j=1}^r$ of $\ker (\calH - \mu)$ so
that $\calB:=\{v\varphi_1 , \ldots , v \varphi_r\}$ is an
orthonormal basis of ${\rm Ran} (S_2)$ (which is
finite-dimensional). This can be done since $v$ is invertible.
Recall that, see (\ref{eq:b0def}),
\[
S_2 b(0) S_2 =
-S_2 v\left[\begin{array}{cc}
\f{- |x-y|}{8\pi}&0\\
0&\f{1}{8\pi} \f{e^{-\sqrt{2\mu}|x-y|}}{\sqrt{2\mu}}
\end{array}
\right] \sigma_3 v S_2.
\]
We denote the matrix of $S_2b(0)S_2$ in the basis $\calB$ by
\[
S_2 b(0) S_2 = \{ a_{ij}\}^r_{i,j=1} =: M.
\]
Then, with
\[
V\varphi_j =: \binom{\alpha_j}{\beta_j},
\]
we have
\begin{align}
a_{ij} &= \langle v \varphi_i , S_2 b(0) S_2 v\varphi_j\rangle \nn\\
&=  \int (\alpha_i(x) , \beta_i(x))
\left[\begin{array}{cc}
\f{|x-y|}{8\pi}&0 \\
0& \f{1}{8\pi} \f{e^{-\sqrt{2\mu}|x-y|}}{\sqrt{2\mu}}
\end{array}
\right]
\binom{\alpha_j(y)}{\beta_j(y)}\;dy\;dx \nn\\
&=  \langle \alpha_i \f{|x-y|}{8\pi}\;, \alpha_j \rangle  +
\lang  \beta_i \f{1}{8\pi} \f{e^{-\sqrt{2\mu}|x-y|}}{\sqrt{2\mu}}\;,\beta_j\rang \nn\\
&= - \langle (-\Delta)^{-1} \alpha_i \;, (-\Delta)^{-1} \alpha_j
\rangle + \langle (-\Delta + 2\mu )^{-1} \beta_i\;,
\left(-\Delta + 2\mu\right)^{-1}\beta_j)\rangle\nn \\
&= -  \lang  \sigma_3 R_0 (\mu) \binom{\alpha_i}{\beta_i}\;,
R_0 (\mu)\binom{\alpha_j}{\beta_j} \rang = -\lang \sigma_3 R_0(\mu)V\varphi_i,R_0(\mu)V\varphi_j\rang\nn \\
&=-\langle \sigma_3 \varphi_i, \varphi_j\rangle. \label{eq:aij}
\end{align}
Here we used that
\[
R_0 (\mu) = \left[\begin{array}{cc}
-\Delta&0\\
0&\Delta-2\mu
\end{array}\right]^{-1} =
\left[\begin{array}{cc}
(-\Delta)^{-1}&0\\
0&(\Delta-2\mu)^{-1}
\end{array}
\right].
\]
Next, write $S_2vR_0(\mu)$ relative to the orthonormal basis
$\calB$:
\begin{align*}
S_2 v R_0(\mu) f &= \sum_j \langle v R_0(\mu) f , v\varphi_j \rangle v \varphi_j\\
&=-\sum_j \langle f , \sigma_3 R_0(\mu) V \varphi_j \rangle v\varphi_j\\
&= \sum_j \langle f , \sigma_3 \varphi_j \rangle v \varphi_j.
\end{align*}
Hence,
\begin{align}
P_\mu f &= -R_0(\mu) \sigma_3 v S_2 (S_2 b(0) S_2)^{-1} S_2 v R_0(\mu) f\nn \\
&=-R_0(\mu) \sigma_3 v S_2 \sum_{i,j} v\varphi_i M^{-1}_{ij}
\langle f , \sigma_3 \varphi_j\rangle\nn \\
&=\sum_{i,j} R_0(\mu) V \varphi_i M^{-1}_{ij} \langle f, \sigma_3 \varphi_j\rangle \nn\\
&= - \sum_{i,j} \varphi_i M^{-1}_{ij} \langle f, \sigma_3 \varphi_j\rangle. \label{eq:B-2}
\end{align}
The following properties hold:
\begin{itemize}
\item[$i$)]$\ran P_\mu \subseteq \ker (\calH - \mu ) = \text{span}
\{\varphi_1,\ldots , \varphi_r\}$
\item[$ii$)]$P_\mu \varphi_k = \varphi_k$
\item[$iii$)]$\ker P_\mu = \ker (\calH^\ast-\mu)^\perp = \text{span}
\{\sigma_3 \varphi_1, \ldots , \sigma_3 \varphi_r \}^\perp$
\end{itemize}
Property $i)$ is immediate from \eqref{eq:B-2}, whereas
 $ii)$ follows from \eqref{eq:aij} and~\eqref{eq:B-2}:
\begin{align*}
P_\mu \varphi_k &= -\sum_{i,j} \varphi_i M^{-1}_{ij}
\langle \varphi_k , \sigma_3 \varphi_j\rangle = \sum_{i,j}
\varphi_i M^{-1}_{ij} a_{j,k}\\
&=\sum_i \varphi_i \delta_{i,k} = \varphi_k.
\end{align*}
Finally, property
$iii)$ can be seen as follows:
\begin{gather*}
P_\mu f  = 0 \Leftrightarrow \langle f , \sigma_3 \varphi_j \rangle = 0 \
\text{for each $j$}\\
\Leftrightarrow f \in \text{span} \{\sigma_3 \varphi_1 , \ldots , \sigma_3 \varphi_r \}^\perp.
\end{gather*}
Since
$ \calH^\ast = \sigma_3 \calH \sigma_3$, we see that
\[ \ker (\calH^\ast - \mu) =
\text{span} \{ \sigma_3 \varphi_1 , \ldots , \sigma_3 \varphi_r \}. \]
The lemma follows.
\end{proof}

Analogously, one obtains expansions around $-\mu$ which involve $P_{-\mu}$.
The previous proposition proves that there is a
direct -- but not orthogonal -- sum representation
\[
L^2 \times L^2 = \ker (\calH -\mu ) \dot{+}
\left[\ker(\calH^\ast - \mu)\right]^\perp
\]
as well as
\[
L^2 \times L^2 = \ker (\calH + \mu ) \dot{+}
\left[\ker(\calH^\ast + \mu)\right]^\perp.
\]
Similarly, for any point $z$ in the discrete spectrum
\[
L^2 \times L^2 = \ker (\calH - z ) \dot{+}
\left[\ker(\calH^\ast - z)\right]^\perp.
\]
Finally, it is a simple matter to check the following:

\begin{lemma}
 The pair-wise products of $P_d, P_\mu$ and $P_{-\mu}$ vanish
where $P_d$ is the Riesz projection
\[ P_d = -\frac{1}{2\pi i}\oint_{\gamma} (\Hil-z)^{-1}\, dz \]
with a simple closed contour $\gamma$ surrounding the discrete
spectrum.
\end{lemma}
\begin{proof}
Suppose that $\Hil f=\mu f$. Then
\[ P_d f = \frac{1}{2\pi i}\oint_{\gamma} (z-\mu)^{-1}\, dz f =0
\]
Hence $P_d P_\mu=0$. Next, suppose that $\Hil f = z f$ and
$\Hil^*g=\mu g$ where $z\in\C$ belongs to the discrete spectrum of
$\Hil$. Then
\[ z\la f,g\ra = \la \Hil f,g\ra = \la f,\Hil^* g\ra = \mu \la
f,g\ra
\]
which implies that $\la f,g\ra=0$. Consequently,
\[ \Ran P_d \subset [\ker (\Hil^*-\mu)]^\perp = \ker P_\mu\]
and thus $P_\mu P_d=0$. The same argument also shows that
\[P_\mu
P_{-\mu}= P_{-\mu}P_\mu = P_dP_{-\mu}=P_{-\mu}P_d =0
\]
and we are done.
\end{proof}

\section{A spectral representation of the evolution and $L^2$
bounds}

The following lemma develops a representation of the (non-unitary)
flow $e^{it\Hil}$ via the resolvents. It relies heavily on the
limiting absorption principle from Section~\ref{sec:gen}. The
statement is of course analogous to the scalar Schr\"odinger case
in which it is a consequence of the spectral theorem and
asymptotic completeness.

\begin{lemma}
\label{lem:rep} Under our assumptions A1)-A4) with $\beta>5$ there
is the representation
\begin{align}
e^{it\Hil} &= \frac{1}{2\pi i}\int_{|\lambda|\ge\mu} e^{it\lambda}\, [(\Hil-(\lambda+i0))^{-1}-(\Hil-(\lambda-i0))^{-1}]\,d\lambda \nn \\
& + \sum_{j} e^{it\Hil} P_{\zeta_j} + e^{it\mu}P_\mu+e^{-it\mu}P_{-\mu},
\label{eq:ac}
\end{align}
where the sum runs over the entire discrete spectrum $\{\zeta_j\}_j$  and
 $P_{\zeta_j}$ is the Riesz projection corresponding to the eigenvalue $\zeta_j$, whereas $P_{\pm\mu}$ are as above.
The formula~\eqref{eq:ac} and the convergence of the integral are
to be understood in the following weak sense: If $\phi,\psi$ belong to
$[W^{2,2}\times W^{2,2}(\R^3)]\cap X_{1+}$, then
\begin{align*}
\la e^{it\Hil}\phi,\psi\ra
&= \lim_{R\to\infty} \frac{1}{2\pi i}\int_{R\ge|\lambda|\ge \mu}\!\! e^{it\lambda}
\big\la [(\Hil-(\lambda+i0))^{-1}-(\Hil-(\lambda-i0))^{-1}]\phi,\psi\big\ra\,d\lambda \\
& \quad + \sum_{j} \la e^{it\Hil}P_{\zeta_j}\phi,\psi \ra
 + e^{it\mu}\la P_\mu\phi,\psi\ra +e^{-it\mu}\la P_{-\mu}\phi,\psi\ra
\end{align*}
for all $t$, where the integrand is well-defined by the limiting absorption principle.
\end{lemma}
\begin{proof}
The evolution $e^{it\Hil}$ is defined via the Hille-Yoshida theorem. Indeed, let $a>0$ be large.
Then $i\Hil-a$ satisfies (with $\rho$ the resolvent set)
\[
 \rho(i\Hil-a)\supset (0,\infty)\text{\ \ and\ \ }\|(i\Hil-a-\lambda)^{-1}\|\le |\lambda|^{-1}\text{\ \ for all\ \ }\lambda>0.
\]
The  inequality can be seen as follows:
\begin{align*}
\|(i\calH-(a+\lambda))^{-1} \| &= \|(i\calH_0-(a+\lambda))^{-1}(I+iV(i\calH_0-(a+\lambda))^{-1})^{-1}\| \\
&\le \frac{1}{a+\lambda} \sum_{k=0}^\infty \Big(\frac{C}{a+\lambda}\Big)^k \le \frac{1}{a-C+\lambda}\le \frac{1}{\lambda},
\end{align*}
as claimed.
Hence $\{e^{t(i\Hil-a)}\}_{t\ge0}$ is a contractive semigroup, so that $\|e^{it\Hil}\|_{2\to2}\le e^{|t|a}$
for all $t\in\R$.
If $\Re z>a$, then there is the Laplace transform
\begin{equation}
\label{eq:Lap}
 (i\Hil-z)^{-1}=-\int_0^\infty e^{-tz}\,e^{it\Hil}\,dt
\end{equation}
as well as its inverse (with $b>a$ and $t>0$)
\begin{equation}
\label{eq:invLap}
e^{it\Hil} = -\frac{1}{2\pi i} \int_{b-i\infty}^{b+i\infty} e^{tz}\, (i\Hil-z)^{-1}\, dz.
\end{equation}
While \eqref{eq:Lap} converges in the norm sense, defining~\eqref{eq:invLap} requires more care.
The claim is that for any $\phi,\psi\in \Dom(\Hil)=W^{2,2}\times W^{2,2}$,
\begin{equation}
\label{eq:limR}
 \la e^{it\Hil} \phi,\psi \ra =
- \lim_{R\to\infty} \frac{1}{2\pi i} \int_{b-iR}^{b+iR} e^{tz}\, \la (i\Hil-z)^{-1}\phi,\psi\ra\, dz.
\end{equation}
To verify this, let $t>0$ and use \eqref{eq:Lap} to conclude that
\bea
-\frac{1}{2\pi i} \int_{b-iR}^{b+iR} e^{tz}\, \la (i\Hil-z)^{-1}\phi,\psi\ra\, dz &=&
 \frac{1}{2\pi i} \int_{b-iR}^{b+iR} e^{tz}\,  \int_0^\infty e^{-sz}\,\la e^{is\Hil}\,\phi,\psi\ra\, dsdz \nn \\
&=& \frac{1}{\pi} \int_0^\infty  e^{(t-s)b}\,\frac{\sin((t-s)R)}{t-s}\,\la e^{is\Hil}\,\phi,\psi\ra\, ds.
\label{eq:dirker}
\eea
Since $e^{(t-s)b}\,\la e^{is\Hil}\,\phi,\psi\ra$ is a $C^1$ function in~$s$ (recall $\phi\in \Dom(\Hil)$)
 as well as  exponentially decaying in $s$ (because of $b>a$),
it follows from standard properties of the Dirichlet kernel that
the limit in~\eqref{eq:dirker} exists and equals $\la e^{it\Hil}\phi,\psi\ra$, as claimed.
Note that if $t<0$, then the limit is zero. Therefore, it follows that for any $b>a$,
\begin{align}
& \la e^{it\Hil} \phi,\psi \ra  \nn \\
&=
- \lim_{R\to\infty} \Big\{
\frac{1}{2\pi i} \int_{b-iR}^{b+iR} e^{tz}\, \la (i\Hil-z)^{-1}\phi,\psi\ra\, dz
- \frac{1}{2\pi i} \int_{-b-iR}^{-b+iR} e^{tz}\, \la (i\Hil-z)^{-1}\phi,\psi\ra\, dz \Big\} \label{eq:2curv}
\end{align}
\begin{figure}
\label{fig:1} \centerline{\hbox{\vbox{ \epsfxsize= 11.0 truecm
\epsfysize= 8.5 truecm \epsfbox{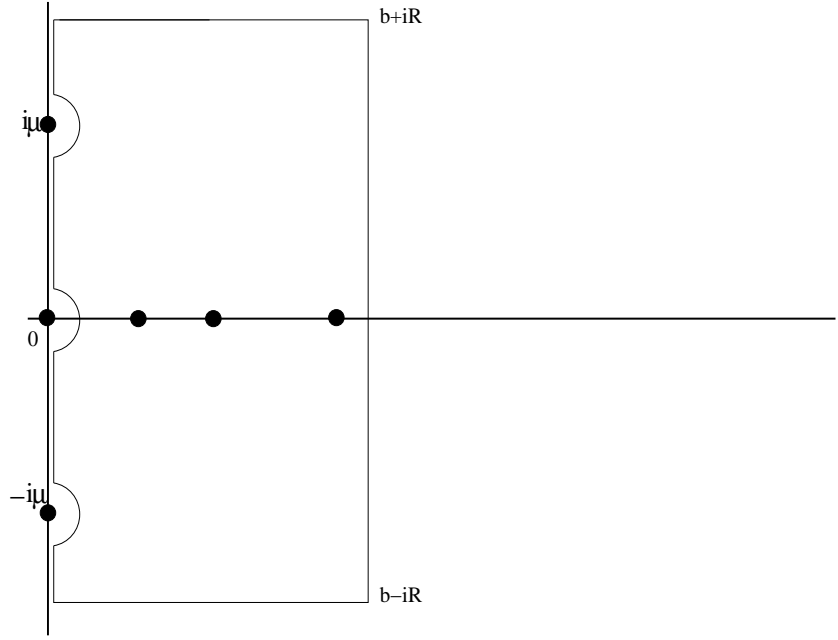}}}}
\caption{The contour $\Gamma^+_{R,\delta}$}
\end{figure}
Consider the contour $\Gamma^+_{R,\delta}$ which is depicted in
Figure~\ref{fig:1}. It has the segment $b-iR$ to $b+iR$ as its
right boundary, and the left boundary contains semi-circular arcs
of radius $\delta>0$ centered at each imaginary
eigenvalue\footnote{This figure depicts the contour for the case
where zero is the only point on the imaginary axis that belongs to
the discrete spectrum of $i\calH$. Generally speaking,
semi-circles surround each point of the discrete spectrum on the
imaginary axis, as well as the edges $\pm i\mu$.} as well as two
semi-circles centered at $\pm i\mu$. Otherwise, the left boundary
abuts on the imaginary axis. Now fix $R$ and some small $\delta>0$
and conclude from the Cauchy theorem that
\[
\frac{1}{2\pi i} \oint_{\Gamma^+_{R,\delta}} e^{tz}\, \la (i\Hil-z)^{-1}\phi,\psi\ra\, dz =
\sum_{j}\frac{1}{2\pi i} \oint_{\gamma_j} e^{tz}\, \la (i\Hil-z)^{-1}\phi,\psi\ra\, dz,
\]
where $\gamma_j$ are small circles (say, of radius $\delta$) around the positive
eigenvalues $\{\lambda_j\}_j$ of $i\Hil$
(in the figure these are indicated by the three dots on the real axis).
Recall that the Riesz projection
\[ P_{\lambda_j} =\frac{1}{2\pi i}\oint_{\gamma_j} (i\Hil-z)^{-1}\, dz \]
satisfies
\[ \Ran (P_{\lambda_j}) = \bigcup_{m=1}^\infty \ker \big[(i\Hil-\lambda_j)^{m}\big] \]
and that the right-hand side stabilizes at some finite (minimal) $M_j=M(\lambda_j)$. I.e.,
\[ \Ran (P_{\lambda_j}) =  \ker \big[(i\Hil-\lambda_j)^{M_j}\big]. \]
This is also the minimal $M_j$ with the property that\footnote{Under our positivity assumption it follows that
the only eigenvalue for which $M_j>1$ is $\lambda=0$. Nevertheless, we still present this argument in general,
since we also want it to apply to $\lambda=0$.}
\begin{equation}
\label{eq:Mest}
 \|(i\Hil-(z-\lambda_j))^{-1}\|\le C\,|z-\lambda_j|^{-M_j}
\end{equation}
as $z\to\lambda_j$.
Now let $p_j(w)$ be the Taylor polynomial of $e^w$ of degree $j$, i.e.,
\[ |e^w-p_j(w)|\le C|w|^{j+1}. \]
Then
\begin{align}
  \frac{1}{2\pi i} \oint_{\gamma_j} e^{tz}\,(i\Hil-z)^{-1}\, dz &=
\frac{1}{2\pi i} \oint_{\gamma_j} e^{t\lambda_j}\big(e^{t(z-\lambda_j)}-p_j(t(z-\lambda_j))\big)(i\Hil-z)^{-1}\, dz \label{eq:cont1}\\
& \quad + \frac{1}{2\pi i} \oint_{\gamma_j} e^{t\lambda_j}p_j(t(z-i\Hil+i\Hil-\lambda_j)) (i\Hil-z)^{-1}\, dz\nn \\
& = \frac{1}{2\pi i} \oint_{\gamma_j} e^{t\lambda_j}p_j(t(i\Hil-\lambda_j)) (i\Hil-z)^{-1}\, dz\nn \\
&= e^{it\Hil} \frac{1}{2\pi i} \oint_{\gamma_j}  (i\Hil-z)^{-1}\, dz = e^{it\Hil}P_{\lambda_j}.\nn
\end{align}
The integral on the right-hand side of \eqref{eq:cont1} is zero by~\eqref{eq:Mest}.
In conclusion,  if we let $\Gamma^-_{R,\delta}$ be the reflection of $\Gamma^+_{R,\delta}$ about the imaginary axis $i\R$, then
\begin{align*}
&\frac{1}{2\pi i} \oint_{\Gamma^+_{R,\delta}} e^{tz}\, (i\Hil-z)^{-1}\, dz
+ \frac{1}{2\pi i} \oint_{\Gamma^-_{R,\delta}} e^{tz}\,(i\Hil-z)^{-1}\, dz \\
&= \sum_{j:\lambda_j\ne0}\frac{1}{2\pi i} \oint_{\gamma_j} e^{tz}\, (i\Hil-z)^{-1}\, dz =
 \sum_{j:\lambda_j\ne0} e^{it\Hil}P_{\lambda_j}.
\end{align*}
Note that by the limiting absorption estimate
\[
\lim_{R\to\infty}\frac{1}{2\pi i} \int_{[iR,b+iR]} e^{tz}\, \la (i\Hil-z)^{-1}\phi,\psi\ra\, dz=0, \]
as well as
\[
\lim_{R\to\infty}\frac{1}{2\pi i} \int_{[-iR,b-iR]} e^{tz}\, \la (i\Hil-z)^{-1}\phi,\psi\ra\, dz=0. \]
Now let $c^+_\delta$  be a contour that is given as follows: Take a straight line $is+\eps$ with $\mu+\delta\le s <\infty$,
then make a circular loop of radius $\delta$ centered at $i\mu$, followed by a straight line $is-\eps$ with the same $s$
as before. Now pass to the limit $\eps\to0$. Similarly with $c^-_\delta$.
Hence, in view of \eqref{eq:2curv} and the preceding,
\begin{align}
\la e^{it\Hil} \phi,\psi \ra  &= \frac{1}{2\pi i}\int_{c^+_\delta} e^{tz} \la (i\Hil-z)^{-1}\phi,\psi\ra \, dz
+ \frac{1}{2\pi i}\int_{c^-_\delta} e^{tz} \la (i\Hil-z)^{-1}\phi,\psi\ra \, dz \nn \\
& \quad + \sum_{j} \la e^{it\Hil}P_{\zeta_j}\phi,\psi\ra \nn \\
& = \frac{1}{2\pi i}\int_{|\lambda|\ge \mu+\delta} e^{it\lambda} [\la (\Hil-(\lambda+i0))^{-1}\phi,\psi\ra-
\la (\Hil-(\lambda-i0))^{-1}\phi,\psi\ra] \, d\lambda \nn \\
& \quad + \sum_{j} \la e^{it\Hil}P_{\zeta_j}\phi,\psi\ra + e^{it\mu}\la P_{\mu}\phi,\psi\ra
+e^{-it\mu}\la P_{-\mu}\phi,\psi\ra+ O(\sqrt{\delta})
\label{eq:sum_rep}
\end{align}
where the sum extends over the entire discrete spectrum $\{\zeta_j\}_j$ of $i\Hil$.
The integrals over infinite intervals are to be interpreted in the principal value sense.
To pass to~\eqref{eq:sum_rep}, we use the asymptotic expansion of the resolvent $(i\Hil-z)^{-1}$ around $\pm i\mu$.
Indeed, by \eqref{eq:RV} and Lemma~\ref{lem:P0} the expansion of $(i\Hil-z)^{-1}$ around $z=i\mu$ is of the form
$(i\mu-z)^{-1}P_{\mu}+O(\sqrt{\delta})$. Sending $\delta\to0$ implies the lemma.
\end{proof}

\begin{remark} The sum over $\zeta_j$ in \eqref{eq:ac} takes the form
\[ \sum_{\zeta_j} e^{it\Hil} P_{\zeta_j} = \sum_{\zeta_j\ne0} e^{it\zeta_j} P_{\zeta_j} + \sum_{k=0}^m \frac{(it)^k}{k!}\Hil^k P_0\]
where $m$ is the minimal positive integer with $\ker(\Hil^m)= \ker(\Hil^{m+1})$. This is due to the fact that
all points $\zeta_j$ in the discrete spectrum other than zero have the property that $\ker(\Hil-\zeta_j)=\ker[(\Hil-\zeta_j)^2]$.
This typically fails for $\zeta_j=0$.
\end{remark}

The previous proposition has the following corollary.

\begin{cor}
\label{cor:L2} Under our assumptions A1)-A4) the following
stability bound holds:
\[ \sup_{t\ge0} \big\|e^{-it\Hil}P_s f\big\|_2 \le C\|f\|_2, \]
where $I-P_s$ is the Riesz-projection corresponding to the
discrete spectrum.
\end{cor}
\begin{proof}
Write $P_s=P_s^+ + P_s^{-}$, where $\pm$ refers to the positive
and negative halves of the essential spectrum, respectively. I.e.,
\[ P_s^+ = \frac{1}{2\pi i} \int_{\Gamma} (\Hil-\lambda)^{-1}
\,d\lambda
\]
with the usual "thermometer" shaped contour surrounding
$[\mu,\infty)$. Then it suffices to prove
\[
\sup_{t\ge0} \big\|e^{-it\Hil}P_s^+ f\big\|_2 \le C\|f\|_2.
\]
Mainly for clarity of exposition we divide the proof into three
cases, namely $S_1=0$, $S_1\ne0, S_2=0$, and finally $S_2\ne0$.
These operators refer to those arising in the expansion of the
resolvent around $\lambda=\mu$.  The first case $S_1=0$ is what is
meant by $\mu$ being neither an eigenvalue nor a resonance.

\medskip
{\bf Case 1: $S_1=0$}

\smallskip
This case has been treated before, and is a consequence of Kato
theory. To see this, write
\[
e^{-it\Hil}P_s = \frac{1}{2\pi i}\int_{|\lambda|\ge\mu}
e^{it\lambda}\,
[(\Hil-(\lambda+i0))^{-1}-(\Hil-(\lambda-i0))^{-1}]\,d\lambda
\]
as explained in the previous lemma.  By the symmetric resolvent
identity~\eqref{mresolve}
\begin{align*}
&(\Hil-(\lambda\pm i0))^{-1} =
(\Hil_0-(\lambda\pm i0))^{-1} \nn \\
& \quad - (\Hil_0-(\lambda\pm i0))^{-1}
v_1(I+v_2(\Hil_0-(\lambda\pm i0))^{-1}v_1)^{-1}v_2
(\Hil_0-(\lambda\pm i0))^{-1}. \nn
\end{align*}
It suffices to show
that
\[ \int \Big| \la (\Hil_0-(\lambda\pm i0))^{-1}
v_1(I+v_2(\Hil_0-(\lambda\pm i0))^{-1}v_1)^{-1}v_2
(\Hil_0-(\lambda\pm i0))^{-1}f,g \ra \Big| \, d\lambda \le
C\|f\|_2\|g\|_2
\]
By our assumption on $\pm\mu$,
\[ \sup_{\lambda}\|(I+v_2(\Hil_0-(\lambda\pm i0))^{-1}v_1)^{-1}\|_{2\to2} <
\infty
\]
Hence, it suffices to show that
\[
\int \Big| \la
v_2(\Hil_0-(\lambda\pm i0))^{-1}f, v_1^*(\Hil_0-(\lambda\mp
i0))^{-1} g \ra \Big| \, d\lambda \le C\|f\|_2\|g\|_2.
\]
However, by Kato's smoothing theory
\[
\int \| v_2(\Hil_0-(\lambda\pm i0))^{-1}f\|_2^2 \, d\lambda \le
C\|f\|_2^2
\]
and similarly for $g$.

\medskip
{\bf Case 2: $S_1\ne0, S_2=0$}

\smallskip
As in the scalar (self-adjoint) case, we do not expect that a
resonance can destroy the $L^2$ bound. We need to check again that
\[
\sup_{t\ge0}\Big\|\int_{\lambda\ge\mu} e^{it\lambda}\,
[(\Hil-(\lambda+i0))^{-1}-(\Hil-(\lambda-i0))^{-1}]\,d\lambda
\Big\|_{2\to2} \le C.
\]
By the Kato theory argument in the previous case
\[
\sup_{t\ge0} \Big\|\int_{\lambda\ge\mu'} e^{it\lambda}\,
[(\Hil-(\lambda+i0))^{-1}-(\Hil-(\lambda-i0))^{-1}]\,d\lambda
\Big\|_{2\to2} \le C(\mu')
\]
for any fixed $\mu'>\mu$. It remains to deal with the integral
over a small interval $\mu'\ge \lambda\ge\mu$. Since $S_2=0$ and
$S_1\ne0$,
\begin{equation} \label{eq:mr}
\mr(z) =-\f{1}{z} \mr_0(z)v_1Sv_2 \mr_0(z) +\mr_0(z) -
\mr_0(z)v_1E(z)v_2\mr_0(z)
\end{equation}
with a uniformly $L^2$-bounded $E(z)$ for small $z$, where we have
set $\lambda=\mu+z^2$ and
\begin{align*}
\mr(z)&=(\Hil-(\mu+z^2+i0))^{-1}, \text{ for } z>0,\\
\mr(z)&=(\Hil-(\mu+z^2-i0))^{-1}=\overline{(\Hil-(\mu+z^2+i0))^{-1}},
\text{ for } z<0
\end{align*}
and similarly for $\mr_0$.  The second term on the right-hand side
of~\eqref{eq:mr} is just the free evolution and thus is bounded in
time. Moreover, the third term can be treated by Kato smoothing
since $\|E(z)\|_{2\to2}$ is bounded for small~$z$. We have reduced
ourselves to showing that
\begin{equation}
\label{eq:red1}
\sup_{t\ge0}\Big\|\int e^{itz^2}\chi(z)\, \mr_0(z)v_1 Sv_2 \mr_0(z)\,
dz  \Big\|_{2\to2} \le C.
\end{equation}
Here $\chi$ is a smooth bump function supported in $[-1,1]$.
Now let $M_j$ for $1\le
j\le 2$ be the matrices
\[ M_1=\bm 1&0 \\ 0&0 \endm,
\quad M_2=\bm 0&0 \\ 0&1 \endm
\]
Denote
\[
U(t) := \int e^{itz^2}\chi(z)\, \mr_0(z)v_1 Sv_2 \mr_0(z)\,
dz.
\]
Then $M_1U(t)M_1$ has the kernel
\begin{equation}
\label{eq:K1}
 \int_{\R^7}
\frac{e^{itz^2}e^{iz(|x-u_1|+|y-u_2|)}\chi(z)}{|x-u_1||y-u_2|}
 r_1(u_1)c_1(u_2) \,du_1du_2\, dz
\end{equation}
where $r_1$ and $c_1$ are in $L^1(\R^3)$. We claim that
\begin{equation}
\label{eq:K2}
\sup_{t\ge0}\|M_1U(t)M_1\|_{2\to2} <\infty.
\end{equation}
In view of~\eqref{eq:K1} this reduces to showing that
\[
\sup_{t\ge0, u_1,u_2}\Big|   \int_{\R^7}
\frac{e^{itz^2}e^{iz(|x -u_1|+|y-u_2 |)}\chi(z)}{|x-u_1 ||y-u_2 |} f(x)g(y)\, dxdydz \Big|\le
C\|f\|_2\|g\|_2
\]
for all $f,g\in L^2(\R^3)$. This in turn is the same as
\begin{equation}
\label{eq:K11}
\sup_{t\ge0}\Big|   \int_{\R^7}
\frac{e^{itz^2}e^{iz(|x  |+|y  |)}\chi(z)}{|x  ||y  |} f(x)g(y)\, dxdy dz  \Big|\le
C\|f\|_2\|g\|_2.
\end{equation}
Let
\[
F(r) =\chi_{[0,\infty)}(r) r \int_{S^2}  f(r\omega)  \sigma(d\omega),\,\,\,\,
G(r) =\chi_{[0,\infty)}(r) r \int_{S^2}  g(r\omega)  \sigma(d\omega).
\]
By Cauchy-Schwarz,
\[
\|F\|_2 \les \|f\|_{L^2(\R^3)}, \qquad \|G\|_2 \les
\|g\|_{L^2(\R^3)}.
\]
The following calculation finishes the proof
of (\ref{eq:K2}):
\begin{align*}
\pref{eq:K11}&= \sup_{t\ge0}\Big|   \int_{\R^3}
 e^{itz^2}e^{iz(r+s)}\chi(z)  F(r)G(s) \, drds dz\Big|\\
 &=\sup_{t\ge0}\Big|   \int_{\R }
 e^{itz^2} \chi(z)  \widehat{F}(z)\widehat{G}(z) \,   dz\Big|\\
 &\leq  \|\widehat{F}\|_{L^2(\R)} \|\widehat{G}\|_{L^2(\R)} = \| F \|_{L^2(\R)} \| G \|_{L^2(\R)}\\
 & \leq C\|f\|_2\|g\|_2.
\end{align*}
Another contribution is given by $M_2 U(t) M_2$.
The kernel here takes the form
\[
 \int_{\R^7}
\frac{e^{itz^2}}{|x-u_1||y-u_2|}\chi(z)
 e^{-\sqrt{2\mu+z^2}(|x-u_1|+|y-u_2|)}r_2(u_1)c_2(u_2) \,du_1du_2\, dz
\]
where $r_2$ and $c_2$ are in $L^2(\R^3)$. The uniform $L^2$ bound
here is even easier, since
\[ \Big\| \int_{\R^3} \frac{e^{-\sqrt{2\mu+z^2}|x-u|}}{|x-u|}
f(u)\, du \Big\|_2 \le C(1+|z|)^{-1} \|f\|_2
\]
so that the desired $L^2$ bound follows by putting $L^2$ norms
inside the integral.

Finally, we claim that both $\|M_1U(t)M_2\|_{2\to2}$ and
$\|M_2U(t)M_1\|_{2\to2}$ are bounded in~$t$. Without loss of
generality we consider the former.

We first remark that
\[
\sup_{t\ge0} \Big\|  \int e^{itz^2}\chi(z)\, M_1\mr_0(z)v_1 Sv_2
(\mr_0(z)- \mr_0(0))M_2\, dz \Big\|_{2\to2} < \infty
\]
by Kato smoothing theory. Indeed, it is easy to check that
\[ \|\chi(z)(\mr_0(z)-\mr_0(0))M_2\|_{2\to2} \les \left\{
\begin{array}{cc} |z| & \text{\ if\ }|z|<1 \\
                    |z|^{-1} & \text{\ if\ }|z|>1
\end{array}
\right. \] Thus we regain the $z$ which we lost due to the
singularity of $\frac{1}{z}S$. It remains to show that
\begin{align}
\sup_{t\ge0} \Big\|  \int e^{itz^2}\chi(z)\, M_1\mr_0(z) v_1 Sv_2
\mr_0(0) M_2\, dz \Big\|_{2\to2} \leq C. \nn
\end{align}
This follows from
\begin{align}\label{eq:m1m2s}
\Big|  \int_{\R^4} e^{itz^2}\frac{e^{iz|x|}}{|x|}\chi(z) f(x)   dxdz \Big|  \leq C\|f\|_2.
\end{align}
As before (with $F(r) =\chi_{[0,\infty)}(r) r \int_{S^2}  f(r\omega)  \sigma(d\omega)$),
\begin{align*}\pref{eq:m1m2s} &=\Big|  \int_{\R} e^{itz^2} \chi(z) \widehat{F}(z)  dz \Big| \\
&\leq \|\chi\|_2  \|\widehat{F}\|_2\les \| F \|_2 \les \|f\|_2
\end{align*}
and we are done with Case~2.

\medskip
{\bf Case 3: $S_2\ne0$}

\smallskip

Let $\Gamma_1(z)=(A(z)+S_1)^{-1}$ and
$\Gamma_2(z)=(m(z)+S_2)^{-1}$. We proved in the previous section
that these are analytic functions for small $z$ and, moreover, for
all small $z\ne0$
\begin{align*}
A(z)^{-1} &= \Gamma_1(z) + \frac{1}{z}
\Gamma_1(z)S_1\Gamma_2(z)S_1\Gamma_1(z) \\
&\quad + \frac{1}{z^2} \Gamma_1(z)S_1 \Gamma_2(z) S_2 b(z)^{-1}
S_2 \Gamma_2(z) S_1 \Gamma_1(z)
\end{align*}
see \eqref{eq:Aexp}. As usual, let $\lambda=\mu+z^2$ with $\Im
z>0$. Then by the symmetric resolvent identity
\begin{align}
 R(\lambda) &:=(\calH - \lambda)^{-1} =
R_0(\lambda) - R_0(\lambda) v_1A(z)^{-1} v_2 R_0(\lambda) \nn\\
 & =R_0(\lambda) - R_0(\lambda) v_1 \Gamma_1(z) v_2 R_0(\lambda)
 -\frac{1}{z}R_0(\lambda) v_1 \Gamma_1(z)S_1\Gamma_2(z)S_1\Gamma_1(z) v_2
 R_0(\lambda)\label{eq:old}\\
 & \quad - \frac{1}{z^2} R_0(\lambda) v_1 \Gamma_1(z)S_1 \Gamma_2(z) S_2 b(z)^{-1}
S_2 \Gamma_2(z) S_1 \Gamma_1(z) v_2 R_0(\lambda) \label{eq:new}
\end{align}
provided $z$ is also small. Note that the contributions of the
terms in~\eqref{eq:old} to the $L^2$ operator norm has been dealt
with in Case~2. Therefore, it suffices to deal
with~\eqref{eq:new}. First, set
\[
T(z)= v_1 \Gamma_1(z)S_1 \Gamma_2(z) S_2 b(z)^{-1}
S_2 \Gamma_2(z) S_1 \Gamma_1(z) v_2
\]
and write
\begin{align}
& \frac{1}{z^2} R_0(\lambda) v_1 \Gamma_1(z)S_1 \Gamma_2(z) S_2
b(z)^{-1} S_2 \Gamma_2(z) S_1 \Gamma_1(z) v_2 R_0(\lambda)
\nn\\
& =\frac{1}{z} R_0(\lambda) \frac{T(z)-T(0)}{z} R_0(\lambda)  +
\frac{1}{z^2} R_0(\lambda) T(0) R_0(\lambda) \label{eq:teilung}
\end{align}
The first term in~\eqref{eq:teilung} only has a $z^{-1}$
singularity, and can therefore be treated as in Case~2. In view of
Lemma~\ref{lem:P0},
\begin{align*}
&R_0(\lambda)T(z)R_0(\lambda)\Big|_{z=0} \\
&=R_0(\lambda) v_1 \Gamma_1(z)S_1 \Gamma_2(z) S_2 b(z)^{-1} S_2
\Gamma_2(z) S_1 \Gamma_1(z) v_2 R_0(\lambda)\bigg|_{z=0} = P
\end{align*}
with $P$ being the projection onto the eigenspace at $\mu$. Now
use the resolvent identity again to conclude that (with $\Im z>0$)
\begin{equation}
\label{eq:RTR} R_0(\lambda) T(0) R_0(\lambda) = P-z^2 R_0(\lambda)
P -z^2 PR_0(\lambda) + z^4 R_0(\lambda)PR_0(\lambda)
\end{equation}
Hence, the contribution of the second term in~\eqref{eq:teilung}
to the $L^2$ operator norm of $e^{it\Hil}$ reduces to
understanding the operator norm of
\[
\int e^{itz^2} z^3\chi(z) \mr_0(z)P\mr_0(z) \, dz
\]
with $\mr_0(z)$ as above (the first three terms on the right-hand
side of \eqref{eq:RTR} are straightforward to deal with). We again
need to consider each of the (essentially scalar) operators
\[
\int e^{itz^2} z^3\chi(z) \mr_0(z)M_jPM_k\mr_0(z) \, dz
\]
for $j,k=1,2$ separately.  According to~\eqref{eq:B-2}
\[
P = - \sum_{i,j} \varphi_i M^{-1}_{ij} \langle \,\cdot\,, \sigma_3 \varphi_j\rangle
\]
where $\{\varphi_j\}_{j=1}^r$ is a suitable basis of
$\ker(\Hil-\mu)$ and $M^{-1}_{ij}$ are some matrix coefficients,
see the proof of Lemma~\ref{lem:P0}. Therefore, the case $j=k=2$
is obvious.   The case $j=k=1$ reduces to establishing that for
any $f,g$ which are the first components of functions in
$\ker(\Hil-\mu)$, we have
\[
\sup_{t\ge0} \Bigl\| \int_{\R^7} e^{itz^2} z^3\chi(z)
\frac{e^{iz(|x-x_1|+|y_1-y|)}}{|x-x_1||y_1-y|}
f(x_1)g(y_1)\,dx_1dy_1\,dz\Bigr\|_{2\to2} \le C
\]
or by duality that
\[
\sup_{t\ge0} \Bigl| \int_{\R^{13}} e^{itz^2} z^3\chi(z)
\frac{e^{iz(|x-x_1|+|y_1-y|)}}{|x-x_1||y_1-y|}
f(x_1)g(y_1)\,dx_1dy_1\,dz\,\phi(x)\psi(y)\,dxdy\Bigr| \le C
\|\phi\|_2\|\psi\|_2
\]
for any pair $\phi,\psi$ of Schwartz functions, say. This is the
same as showing that
\[
\sup_{t\ge0} \Bigl| \int_{\R^{7}} e^{itz^2} z^3\chi(z)
\frac{e^{iz(|x|+|y|)}}{|x||y|}
(f\ast\phi)(x)(g\ast\psi)(y)\,dxdy\,dz\Bigr| \le C
\|\phi\|_2\|\psi\|_2
\]
We remark that this estimate is different from the ones we
encountered in Case~2 since $f,g\in L^2$ but not necessarily
$f,g\in L^1$. We set
\[
F(r) = r \chi_{[r>0]}\int_{S^2}
(f\ast\phi)(r\omega)\,\sigma(d\omega),\quad G(r) = r
\chi_{[r>0]}\int_{S^2} (g\ast\psi)(r\omega)\,\sigma(d\omega).
\]
Since $\phi,\psi\in L^1(\R^3)$, we conclude that $F,G \in L^2_r$.
Moreover, $\partial_r F, \partial_r G \in L^2$ and
$\widehat{\partial_r F}(z)=z\hat{F}(z)$, $\widehat{\partial_r
G}(z)=z\hat{G}(z)$. To see that $\partial_r F\in L^2$, observe
that
\[ \|\partial_r F\|_2 \les (\|\nabla (f\ast\phi)\|_2 +
\|\,|x|^{-1}(f\ast\phi)\|_2) \les \|f\ast\nabla\phi\|_2 \les
\|f\|_2 \|\nabla \phi\|_1
\]
where we applied Hardy's inequality in the second step. Hence, we
need to show that
\[
\sup_{t\ge0} \Bigl| \int e^{itz^2} z\chi(z) \widehat{\partial_r
F}(z)  \widehat{\partial_r G}(z)\,dz\Bigr| \le C
\|\phi\|_2\|\psi\|_2
\]
which in turn reduces to proving that
\[
\|\partial_r F\|_2 \le C\|\phi\|_2,\qquad \|\partial_r G\|_2
\le C\|\psi\|_2.
\]
For this it suffices to check that
\[
\|\nabla (f \ast \phi)\|_2 + \||x|^{-1} (f\ast \phi)\|_2 \le
C\|\phi\|_2.
\]
By Hardy's inequality the second term on the
left-hand side is controlled by the first. Hence, we need to show
that
\[
\|\widehat{\nabla f}\|_\infty \le C.
\]
We recall that $f$ is assumed to satisfy $f\in L^2(\R^3)$ and
\[ -\Delta f = V_1 f + V_2 \tilde{f} \]
for some $\tilde{f}\in L^2(\R^3)$. By the assumed decay of $V_1,
V_2$ we have $V_1 f + V_2 \tilde{f}\in L^1(\R^3)$ and also
\[ \int_{\R^3} (V_1 f + V_2 \tilde{f})\, dx=0\]
see Lemma~\ref{L:ms1s2}.  Set $h=V_1 f + V_2 \tilde{f}$. It
follows that $|\hat{h}(\xi)|\le C|\xi|$ and thus $\nabla f= \nabla
G_0 h $ implies that
\[ \sup_{\xi\in\R^3}|\widehat{\nabla f}(\xi)|\le C,\]
as desired.

It remains to consider the case $j=1, k=2$ ($j=2, k=1$ being
symmetric). We need to show that for any Schwartz functions
$\phi,\psi$
\[ \sup_{t\ge0} \Bigl| \int_{\R^{7}} e^{itz^2}
z^3\chi(z) \frac{e^{iz|x|-|y|\sqrt{2\mu+z^2}}}{|x||y|}
(f\ast\phi)(x)(g\ast\psi)(y)\,dxdy\,dz\Bigr| \le C
\|\phi\|_2\|\psi\|_2
\]
where $f,g$ are the first and second components, respectively, of
functions in $\ker(\Hil-\mu)$. Since
\[ \int_{\R^3} \frac{e^{-|y|\sqrt{2\mu+z^2}}}{ |y|}|(g\ast\psi)(y)|\,dy
\les \|g\ast \psi\|_\infty \les \|g\|_2 \|\psi\|_2,\] we see that
we are reduced to showing that
\[
 \int z^2 \chi(z) |\widehat{\partial_r F}(z)|\, dz
 \le C\|\phi\|_2.
\]
This, however, was already established for the case $j=k=1$, and
we are done.
\end{proof}

\section{$L^1\to L^\infty$ estimates}
In this section, we bound the $L^1\to L^\infty$ operator norm of
$e^{it\Hil}P_s=e^{it\Hil}P_s^++e^{it\Hil}P_s^-$. In view of the
previous section, the kernel of $e^{it\Hil}P_s^+$, truncated to
energies close to $\mu$, is
\[
K_t:=\frac{1}{2\pi i}\int_{\lambda>\mu}e^{it\lambda}\chi(\lambda-\mu)
[(\Hil-(\lambda+i0))^{-1}-(\Hil-(\lambda-i0))^{-1}] d\lambda,
\]
where $\chi$ is an even Schwartz function supported in
$(-\lambda_0,\lambda_0)$  and identically equal to $1$ in
$(-\lambda_0/2,\lambda_0/2)$. Here $\lambda_0$ is a small constant
which will be determined later. The dispersive estimates for the
remaining operators, i.e., those defined in terms of $1-\chi$, were
obtained in~\cite{Sch1}. Hence, we shall only work with energies
close to the thresholds $\pm\mu$. By a simple change of variable and
redefining $\chi$, we have \bea \label{Kt} K_t=\frac{1}{\pi i}
e^{it\mu} \int_{-\infty}^\infty e^{it z^2} z \chi(z) \mr(z) dz, \eea
where
\begin{align*}
\mr(z)&=(\Hil-(\mu+z^2+i0))^{-1}, \text{ for } z>0,\\
\mr(z)&=(\Hil-(\mu+z^2-i0))^{-1}=\overline{(\Hil-(\mu+z^2+i0))^{-1}}, \text{ for } z<0.
\end{align*}
We also define
\begin{align*}
\mr_0(z)&=R_0(\mu+z^2), \text{ for } z>0,\\
\mr_0(z)&=\overline{R_0(\mu+z^2)}, \text{ for } z<0.
\end{align*}
Note that, with this definition, for all $z\in \R$ we have
$$
\mr_0(z)(x,y)=\f{1}{4\pi|x-y|}\left[
\begin{array}{cc}
 e^{iz|x-y|} & 0\\
0&- e^{-\sqrt{2\mu+z^2}|x-y|}
\end{array}
\right]
$$
and
\bea \label{resexp}
\mr(z)=\mr_0(z) - \mr_0(z) v_1A(z)^{-1}v_2 \mr_0(z),
\eea
where
$$
A(z)=I+v_2\mr_0(z)v_1.
$$
We will use the following simple lemma repeatedly. It is  used in
\cite{RodSch,goldbergschlag, ErdSch}.
\begin{lemma}\label{L:birb}
Let $F\in L^1$ be differentiable on $\R$ with $F^\prime \in L^1$. Then
\begin{align*}
\text{i) }&\left|\int_{-\infty}^\infty e^{itz^2} F(z) dz\right|\lesssim t^{-1/2} \|\widehat{F}\|_1, \\
\text{ii) }&\left|\int_{-\infty}^\infty e^{itz^2} z  F(z) dz\right|\lesssim t^{-3/2} \|\widehat{F^\prime}\|_1.
\end{align*}
\end{lemma}
\begin{proof}
This follows from
\[ |\widehat{e^{itz^2}}(u)| = c|t|^{-\half}\]
and Parseval's identity.
\end{proof}

\subsection{$\mu$ is a resonance but not an eigenvalue}

Now, we prove Theorem~\ref{T:scalar2} when $\mu$ is a resonance but not an eigenvalue.
In this case $S_2=0$ and we have
\begin{align*}
 A(z)^{-1}  = \left(A(z) + S_1\right)^{-1}  + \f{1}{z} \left(A(z) + S_1\right)^{-1} S_1
 m(z)^{-1} S_1 \left(A(z) + S_1\right)^{-1},
\end{align*}
where
\begin{align}
\left(A(z) + S_1\right)^{-1}
&=(A_0 + S_1)^{-1} + \sum_{k=1}^\infty (-1)^k z^k (A_0 + S_1)^{-1}
\left[A_1(z) (A_0+S_1)^{-1}\right]^k\nn\\
&=:(A_0+S_1)^{-1}+z E_1(z),\nn \\
m(z)^{-1}&=m(0)^{-1}+\sum_{k=1}^\infty (-1)^k z^k  m(0)^{-1}
\left[m_1(z) m(0)^{-1}\right]^k \label{eq:m0inv}\\
&=:m(0)^{-1}+z E_2(z).\nn
\end{align}
Thus, using \pref{s1commute}, we obtain\footnote{We note that the corresponding equation (22)
in \cite{ErdSch} has a couple of misprints. It should be replaced with the equation (\ref{S2=0expan}).
The rest of the proof in \cite{ErdSch} is not affected by this change.}
\begin{align}
\label{S2=0expan}
A(z)^{-1}  & =  \frac{1}{z}S_1  m(0)^{-1} S_1  \\
& + \left(A(z) + S_1\right)^{-1} \nonumber
\\
 &+   E_1(z) S_1m(z)^{-1} S_1 \left(A(z) + S_1\right)^{-1}\nonumber \\
 &+    S_1
 E_2(z) S_1 \left(A(z) + S_1\right)^{-1} \nonumber\\
 &+
  S_1
 m(0)^{-1} S_1 E_1(z)  \nonumber \\
& =: \f{1}{z}S+E(z).\nonumber
\end{align}
Note that $m(0)=-\frac{i}{4\pi} S_1 vP_1 vS_1$ is invertible in
$S_1L^2$. Since $P_1$ is of rank one, both $m(0)$ and $S_1$ are
rank one operators. Let $S_1(x,y)=\varphi(x)\varphi^*(y)$, where
$\varphi$ is the unique function  satisfying i) $\|\varphi\|_2=1$,
ii) $\varphi=v_2g=vg$ for a resonance function $g$, as well as
iii) $P_1v\varphi=\binom{c}{0}$ with $c>0$ (see
Lemma~\ref{L:ms1s2}).

Using $c$  in the definition of $m(0)$, it is easy to see
that\footnote{$\varphi$ is real-valued, so $\varphi^*=\varphi^t$}
\bea \label{SDEF} S(x,y)=\frac{4\pi i}{c^2}S_1(x,y)=\frac{4\pi
i}{c^2}\varphi(x)\varphi^*(y). \eea Plugging (\ref{S2=0expan}) and
(\ref{SDEF}) into (\ref{resexp}), we have
\begin{align}
\mr(z) &=-\frac{4\pi i}{c^2z} \mr_0(z)v_1\varphi\varphi^*v_2 \mr_0(z) \nn\\
&+\mr_0(z) - \mr_0(z)v_1E(z)v_2\mr_0(z) \nn \\
&= \frac{4\pi i}{c^2z} \mr_0(z)\sigma_3 \psi\psi^*  \mr_0(z)\label{eq:S1exp}\\
&+\mr_0(z) - \mr_0(z)v_1E(z)v_2\mr_0(z),\nn
\end{align}
where $\psi=v\varphi=-\sigma_3 Vg \in
L^{2,\beta-\frac{1}{2}-}\subset L^1\cap L^2(\R^3) $. Using this in
(\ref{Kt}), we get
$$
K_t(x,y)= \frac{1}{\pi i} e^{it\mu}(\frac{4\pi i}{c^2} K_1(x,y) +  K_2(x,y) - K_3(x,y)),
$$
where
\begin{align}
K_1(x,y)&  = \int_{-\infty}^\infty   e^{itz^2} \chi(z)
[\mr_0(z)\sigma_3\psi\psi^*\mr_0(z)](x,y) dz,\nn \\
K_2(x,y)& = \int_{-\infty}^\infty e^{itz^2}z \chi(z) \mr_0(z)(x,y)  d z \nn \\
K_3(x,y)& = \int_{-\infty}^\infty e^{itz^2}z \chi(z)[\mr_0(z)v_1E(z)v_2\mr_0(z)](x,y)  d z.
\label{eq:K3def}
\end{align}
First,  we deal with $K_1$. Let $\psi=\left[\begin{array}{c}\psi_1\\\psi_2\end{array}\right]$.
\begin{align}
\label{k1}
K_1(x,y)     = \frac{1}{16\pi^2}\left[\begin{array}{cc}
K_1^{11} & K_1^{12} \\ K_1^{21} & K_1^{22}
\end{array}\right],
\end{align}
where
\begin{align*}
K_1^{11} & = \int_{\R^7} \f{e^{itz^2} \chi(z)}{|x-u_1||y-u_2|}
e^{iz(|x-u_1|+|y-u_2|)}\psi_1(u_1){\psi_1}(u_2) d u_1 d u_2   dz\\
K_1^{12} & =-\int_{\R^7} \f{e^{itz^2} \chi(z)}{|x-u_1||y-u_2|}
e^{iz |x-u_1|-\sqrt{2\mu+z^2}|y-u_2| }\psi_1(u_1){\psi_2}(u_2)  d u_1 d u_2   dz\\
K_1^{21} & =\int_{\R^7} \f{e^{itz^2} \chi(z)}{|x-u_1||y-u_2|}
e^{iz |y-u_2|-\sqrt{2\mu+z^2}|x-u_1|}\psi_2(u_1){\psi_1}(u_2) d u_1 d u_2   dz\\
K_1^{22} & =-\int_{\R^7} \f{e^{itz^2} \chi(z)}{|x-u_1||y-u_2|}
e^{-\sqrt{2\mu+z^2}(|x-u_1|+|y-u_2|)}\psi_2(u_1){\psi_2}(u_2)  d
u_1 d u_2   dz.
\end{align*}

As in the scalar case, we will prove that $K_1$ is a sum of two
operators the first one is of finite rank and its $L^1\rightarrow
L^\infty$ norm decays like $t^{-1/2}$, the second one is
dispersive, i.e., its $L^1\rightarrow L^\infty$ norm decays like
$t^{-3/2}$. It suffices to prove this claim for each of the
components of $K_1$.

 First we consider $K_1^{11}$. Let $a_1=|x-u_1|$, $a_2=|y-u_2|$ and $a=a_1+a_2$.
\begin{align} \nn
 K_1^{11}=\int_{\R^7}  e^{itz^2} \chi(z)\f{\cos(za)}{a_1a_2}\psi_1(u_1){\psi_1}(u_2)  d u_1 d u_2   dz.
\end{align}
We have
(see \cite{ErdSch})
\begin{align} \label{fourier}
 \int_{-\infty}^\infty   e^{itz^2} \chi(z) \cos(z a)   dz
 =&\frac{e^{i\pi/4}}{\sqrt{4\pi t}}\int_{-\infty}^\infty e^{-i(u^2+a^2)/4t}  \widehat{\chi}(u)du \\
&+\frac{e^{i\pi/4}}{\sqrt{4\pi t}}\int_{-\infty}^\infty e^{-i(u^2+a^2)/4t} (\cos(\f{ua}{2t})-1)\widehat{\chi}(u)du\nonumber\\
=
&\frac{e^{i\pi/4}}{\sqrt{4\pi t}}\int_{-\infty}^\infty e^{-i(u^2+a_1^2+a_2^2)/4t}  \widehat{\chi}(u)du\nonumber\\
&+\frac{e^{i\pi/4}}{\sqrt{4\pi t}}\int_{-\infty}^\infty e^{-iu^2/4t} (e^{-ia^2/4t}-e^{-i(a_1^2+a_2^2)/4t})  \widehat{\chi}(u)du\nonumber\\
&+\frac{e^{i\pi/4}}{\sqrt{4\pi t}}\int_{-\infty}^\infty e^{-i(u^2+a^2)/4t} (\cos(\f{ua}{2t})-1)\widehat{\chi}(u)du\nonumber\\
=:&C_1+C_2 +C_3.\nonumber
\end{align}
In \cite{ErdSch}, we proved that the contribution of $C_2$ and $C_3$ in $K_1^{11}$ are dispersive, see pages 367-369 in that paper.
The contribution of $C_1$ is
\begin{align*}
t^{-1/2}F_{11}(t) :=\frac{ h(t)}{\sqrt{t}} \Big[\int_{\R^3}
\f{e^{-i|x-u_1|^2/4t} \psi_1(u_1)}{|x-u_1|}d u_1 \Big]\Big[
\int_{\R^3}\f{e^{-i|y-u_2|^2/4t} {\psi_1}(u_2)}{|y-u_2|}  d
u_2\Big],
\end{align*}
where $h(t)=\frac{e^{i\pi/4} }{\sqrt{4\pi }}\int_{-\infty}^\infty e^{-i u^2/4t}  \widehat{\chi}(u)du$.

Now, we consider $K_1^{22}$, the others can be treated similarly. Let $a=|x-u_1|+|y-u_2|$.
\begin{align*}
K_1^{22}(x,y)
&=-\int_{-\infty}^\infty \int_{\R^6} \f{e^{itz^2} \chi(z)}{|x-u_1||y-u_2|}
  e^{-a\sqrt{2\mu}} \psi_2(u_1){\psi_2}(u_2)
   d u_1 d u_2   dz \\
& - \int_{-\infty}^\infty \int_{\R^6} \f{e^{itz^2} z a \chi(z)}{|x-u_1||y-u_2|}
  \f{e^{-a\sqrt{2\mu+z^2}} - e^{-a\sqrt{2\mu}}}{za}  \psi_2(u_1){\psi_2}(u_2)
   d u_1 d u_2   dz \\
& =: t^{-1/2} F_{22}(x,y) + K_{122}(x,y).
\end{align*}
By Lemma~\ref{L:birb}, $\|F_{22}\|_{1\rightarrow \infty} \lesssim 1$.
Before we prove that $K_{122}$ is dispersive,
we note that the kernel of the operator $F_t$ in Theorem~\ref{T:scalar2}
when $\mu$ is a resonance but not an eigenvalue is
\begin{align}\label{ft_s2=0}
F_t(x,y)&=\frac{e^{it\mu}h(t)}{4\pi^2c^2}\left[
\begin{array}{cc}
T_t(\psi_1)(x)T_t({\psi_1})(y) & -T_t(\psi_1)(x) Q({\psi_2})(y)\\
Q(\psi_2)(x)T_t({\psi_1})(y)&- Q(\psi_2)(x)Q({\psi_2})(y)
\end{array}\right] \\
&=\frac{e^{it\mu}h(t)}{4\pi^2c^2}\left[
\begin{array}{c}
T_t(\psi_1)(x) \\
Q(\psi_2)(x)
\end{array}\right]
\left[ \begin{array}{cc} T_t({\psi_1})(y) &
-Q({\psi_2})(y)\end{array}\right],\nn
\end{align}
where
$$
T_t(f)(x):=\int\frac{e^{-i|x-u|^2/4t}f(u)}{|x-u|}du,\,\,\,\,\,\,\,Q(f)(x):=\int\frac{e^{-\sqrt{2\mu}|x-u|}f(u)}{|x-u|}du.
$$

To prove that $K_{122}$ is dispersive we need the following
calculus lemma.

\begin{lemma}
\label{lem:calc} For any $k\in\R$ define
\[
g_k(x) = \frac{x}{\sqrt{k^2+x^2}} e^{-\sqrt{x^2+k^2}}.
\]
Then
\[
\|g_k\|_1 + \|g_k'\|_1 + |k| \|g_k''\|_1 + |k|^2 \|g_k'''\|_1 \les P(k) e^{-|k|},
\]
where $P$ is a polynomial in $k$.
\end{lemma}
\begin{proof} Clearly,
\[
g_k(x) = - \frac{d}{dx} e^{-\sqrt{x^2+k^2}}.
\]
Hence
\[
\|g_k\|_1 = 2\int_0^\infty g_k(x)\, dx = 2e^{-|k|}.
\]
Next, $g_k(x)=\widetilde{g_k}(x/k)$ where
\[
\widetilde{g_k}(x)=\frac{x}{\sqrt{1+x^2}}e^{-|k|\sqrt{1+x^2}}.
\]
Hence, \bea \label{gktilde} \|g_k^{(j)}\|_1 = |k|^{-(j-1)}
\|\widetilde{g_k}^{(j)}\|_1, \,\,\,\, j=1,2,... \eea Note that all
derivatives of $\sqrt{1+x^2}$ are bounded functions. Therefore, by
Leibnitz's formula
\[
|\widetilde{g_k}^{(j)}(x)| \lesssim \Big|\big(\frac{x}{\sqrt{1+x^2}}\big)^{(j)}\Big| e^{-|k|}
+ |k|(1+|k|^{(j-1)}) e^{-|k|\sqrt{1+x^2}}.
\]
Note that all derivatives of $\frac{x}{\sqrt{1+x^2}}$ are in $L^1$.
Thus,
\[
\|\widetilde{g_k}^{(j)}\|_1\lesssim e^{-|k|}+(1+|k|^{(j-1)}) \||k| e^{-|k|\sqrt{1+x^2}}\|_1
\lesssim P_j(k) e^{-|k|},
\]
where $P_j(k)$ is a polynomial.
Using this in (\ref{gktilde}) yields the assertion of the lemma.
\end{proof}

Let $h_a(z):=\f{z}{\sqrt{2\mu a^2+z^2}}e^{-\sqrt{2\mu a^2+z^2}}$.
In view of Lemma~\ref{lem:calc} we have the following bounds for
$h_a$:
\begin{align*}
\|h_a\|_1+\|h_a^\prime\|_1 & \lesssim e^{-|a|\sqrt{\mu}},   \\
\|h_a^{\prime\prime}\|_1  & \lesssim \f{e^{-|a|\sqrt{\mu}}}{\sqrt{\mu}|a|},  \\
\|h_a^{\prime\prime\prime}\|_1  & \lesssim
\f{e^{-|a|\sqrt{\mu}}}{\mu|a|^2}.
\end{align*}
Therefore, we have
\begin{align}
|\widehat{h_a}(\eta)|  \lesssim e^{-\mu |a|} \min(\f{1}{\la\eta\ra}, \f{1}{|a|\la\eta\ra^2}, \f{1}{a^2 \la\eta\ra^3}). \label{hb}
\end{align}
Note that
\begin{align*}
 K_{122}(x,y) &=  - \int_{-\infty}^\infty \int_{\R^6} \int_0^1 \f{e^{itz^2} z a \chi(z)}{|x-u_1||y-u_2|}  h_a(abz)
   \psi_2(u_1){\psi_2}(u_2)
  db d u_1 d u_2   dz.
\end{align*}
Lemma~\ref{L:birb} implies that $ K_{122}$ is dispersive  if we can prove that
\begin{align} \label{K1hat}
 \sup_{x,y}\int_{\R^6} \int_0^1  \|\widehat{(\chi h_a(ab(\cdot)))^\prime}\|_1
   \f{a \psi_2(u_1){\psi_2}(u_2)}{|x-u_1||y-u_2|}
  db d u_1 d u_2
\end{align}
is finite.
Using the Schwartz decay of $\widehat{\chi}$ and (\ref{hb}), we obtain
\begin{align*}
\|\widehat{(\chi h_a(ab(\cdot)))^\prime}\|_1& = 2 \pi \int  |\xi|  \left|\int \widehat\chi(\xi-ba\eta)\widehat{h_a}(\eta)d\eta\right|d \xi \\
&\lesssim  \int   (1+ |ba   \eta|) |\widehat{h_a}(\eta)|d\eta \\
&\lesssim \f{1}{a } e^{-\sqrt{\mu}a}.
\end{align*}
Using this in (\ref{K1hat}), we have
\begin{align*} \pref{K1hat} \lesssim
 \sup_{x,y}\int_{\R^6}  e^{-\sqrt{\mu}a}
   \f{ \psi_2(u_1){\psi_2}(u_2)}{  |x-u_1||y-u_2|}
  d u_1 d u_2  <\infty.
\end{align*}
This finishes the analysis of $K_1$. Note that $K_2$ is the low
energy part of the free evolution and hence it is dispersive. Now,
we consider $K_3$. Let $I_1$ be the first coordinate projection
with the matrix $I_1=[1\;\;0]$ and $I_2=[0\;\;1]$ the second
coordinate projection. We have
\begin{align*}
K_3(x,y) = \f{1}{16\pi^2} \left[\begin{array}{cc}
K_3^{11}(x,y)&- K_3^{12}(x,y)\\
- K_3^{21}(x,y) & K_3^{22}(x,y)
\end{array}\right],
\end{align*}
where
\begin{align*}
K_3^{11}(x,y)&=\int_{-\infty}^\infty \int_{\R^6} \f{e^{itz^2} z \chi(z)e^{iz(|x-u_1|+|y-u_2|)}}{|x-u_1||y-u_2|}
(I_1v_1E(z)v_2I_1^T)(u_1,u_2) d u_1 d u_2   dz \\
K_3^{12}(x,y)&=\int_{-\infty}^\infty \int_{\R^6} \f{e^{itz^2} z \chi(z)e^{iz |x-u_1|-\sqrt{2\mu+z^2}|y-u_2| }}{|x-u_1||y-u_2|}
(I_1v_1E(z)v_2I_2^T)(u_1,u_2) d u_1 d u_2   dz \\
K_3^{21}(x,y)&=\int_{-\infty}^\infty \int_{\R^6} \f{e^{itz^2} z \chi(z)e^{iz |y-u_2|-\sqrt{2\mu+z^2}|x-u_1| }}{|x-u_1||y-u_2|}
(I_2v_1E(z)v_2I_1^T)(u_1,u_2) d u_1 d u_2   dz \\
K_3^{22}(x,y)&=\int_{-\infty}^\infty \int_{\R^6} \f{e^{itz^2} z \chi(z)e^{-\sqrt{2\mu+z^2}(|x-u_1|+|y-u_2|)}}{|x-u_1||y-u_2|}
(I_2v_1E(z)v_2I_2^T)(u_1,u_2) d u_1 d u_2   dz
\end{align*}

The rest of this section is devoted to the proof of
\begin{equation}
\label{eq:K3est}
\sup_{x,y}|K_3^{ij}(x,y)|\lesssim t^{-3/2},\,\, i,j=1,2.
\end{equation}
First, we consider $K_3^{11}$.
Denote
$$\int_{\R^6}\f{d}{dz}\left(  \chi(z)e^{iz(|x-u_1|+|y-u_2|)}
(I_1v_1E(z)v_2I_1^T)(u_1,u_2)\right) \f{ d u_1 d u_2}{|x-u_1||y-u_2|}  $$
by $\calF_{x,y}(z)$.
By Lemma~\ref{L:birb}, it suffices to prove that
\bea\label{fl1}
\sup_{x,y} \| \widehat{\calF_{x,y}}\|_{L^1}<\infty.
\eea

Let us  concentrate on the term where the derivative hits $\chi (I_1v_1E v_2I_1^T)$ (the term
where the derivative hits the exponential is similar):
$$
\tilde{\calF}_{x,y}(z)=\int_{\R^6}[\chi(z)(I_1v_1E(z)v_2I_1^T)]^\prime(u_1,u_2)
\f{e^{iz(|x-u_1|+|y-u_2|)}}{|x-u_1||y-u_2|}du_1du_2.
$$
Note that
\begin{align}
\|\widehat{\tilde{\calF}_{x,y}}(\xi)\|_{L^1}&=\int_{-\infty}^\infty \left| \int_{\R^6}
\widehat{ [\chi(I_1v_1Ev_2I_1^T)]^\prime }(\xi-|x-u_1|-|y-u_2|)(u_1,u_2)
\f{du_1du_2 }{|x-u_1||y-u_2|}\right|d\xi\nn\\
&\leq \int_{\R^6} \int_{-\infty}^\infty \left|\widehat{ [\chi(I_1v_1Ev_2I_1^T)]^\prime }(\xi-|x-u_1|-|y-u_2|)(u_1,u_2)\right|
\f{ d\xi du_1du_2 }{|x-u_1||y-u_2|} \nn\\
&=\int_{\R^6} \int_{-\infty}^\infty \left|\widehat{ [\chi(I_1v_1Ev_2I_1^T)]^\prime }(\xi)(u_1,u_2)\right|
\f{ d\xi du_1du_2}{|x-u_1||y-u_2|}  \label{Fvv}
\end{align}
The second line follows from Minkowski's inequality and Fubini's theorem, the third line follows from a change of variable.
Note that
$I_1v_1E(z)v_2I_1^T(u_1,u_2)$ is a sum of kernels of the form
\[
w_1(u_1)E_{ij}(z)(u_1,u_2)w_2(u_2),\,\,i,j=1,2,\,\,w_1,w_2\in L^2(\R^3)\cap L^\infty(\R^3).
\]
Using this and the inequality (for $w\in L^2(\R^3)\cap L^\infty(\R^3)$)
\begin{align}
\label{vcalc}
\left\| \f{|w(\cdot)| }{|x-\cdot| }\right\|_2^2& =\int_{|x-u|<1} \f{|w(u)|^2 }{|x-u|^2 } du + \int_{|x-u|>1} \f{|w(u)|^2 }{|x-u|^2 } du
\\
&\lesssim \int_{|u|<1} \f{1}{|u|^2 } du + \int_{\R^3}  |w(u)|^2   du \lesssim 1,
\nonumber
\end{align}
in (\ref{Fvv}), we have
\begin{align*}
\|\widehat{\tilde{\calF}_{x,y}}(\xi)\|_{L^1}& \leq \sum_{i,j=1}^2
\int_{-\infty}^\infty \left\|\left|\widehat{(\chi E_{ij} )^\prime}(\xi) \right| \right\|_{L^2\rightarrow L^2}d\xi.
\end{align*}
Therefore, for $\tilde{\calF}_{x,y}$, (\ref{fl1}) follows from
\beeq
\int_{-\infty}^\infty
\left\|\,\left|\widehat{(\chi E_{ij} )^\prime}(\xi) \right|\, \right\|_{L^2\rightarrow L^2}d\xi
< \infty,\,\,\,i,j=1,2.
\label{eq:K3end}
\eneq
We shall use the following elementary lemma from \cite{ErdSch}.

\begin{lemma}\label{L:convolve}
For each $z\in \R$, let $F_1(z)$ and $F_2(z)$ be bounded operators from $L^2(\R^3)$ to $L^2(\R^3)$
with kernels $K_1(z)$ and $K_2(z)$. Suppose that $K_1, K_2$ both have compact support in $z$ and
that $K_j(\cdot)(x,y)\in L^1(\R)$ for a.e.~$x,y\in\R^3$.
Let $F(z)=F_1(z)\circ F_2(z)$ with kernel $K(z)$. Then
$$
\int_{-\infty}^\infty \left\|\,\left|\widehat{K}(\xi)\right|\,\right\|_{ 2\rightarrow 2} d\xi
\leq \left[\int_{-\infty}^\infty \left\|\, \left|\widehat{K_1 }(\xi)\right|\,\right\|_{ 2\rightarrow 2}d\xi\right]
\left[\int_{-\infty}^\infty \left\|\, \left|\widehat{K_2 }(\xi)\right|\,\right\|_{ 2\rightarrow 2}d\xi\right].
$$
\end{lemma}
\noindent

Note that $\f{d}{dz}[\chi(z)E_{ij}(z)]$ is a sum of operators each
of which is a composition of  operators from the list below  (here
$\chi$ is a suitably chosen smooth cut-off supported in a small
neighborhood of the origin):
\begin{align*}
F_1(z) &=  \chi(z) (A(z)+S_1)^{-1}, \\
F_2(z) &=  \chi(z) E_1(z),\\
F_3(z) &=  \chi(z) S_1 m(z)^{-1}S_1, \\
F_4(z) &=  \chi(z) S_1 E_2(z) S_1,
\end{align*}
and their $z$ derivatives and appropriate projections.
Moreover, we leave it to the reader to check that for each of the combinations that
contribute to $E_{ij}(z)$ the hypotheses of Lemma~\ref{L:convolve} are fulfilled.
Therefore, in light of Lemma~\ref{L:convolve}, the following lemma completes the analysis of $K_3^{11}$.
\begin{lemma}\label{L:hammal}
For each of the operators $F_j$, $j=1,2,3,4$ above,
\bea \label{hammal}
\int_{-\infty}^\infty \left\| \left|\widehat{F_j }(\xi)\right|\right\|_{ 2\rightarrow 2}d\xi < \infty.
\eea
The same statement is valid for their $z$ derivatives, too.
\end{lemma}
\begin{proof}
We omit the analysis of $F_1$ and $F_3$. Recall that
\begin{align*}
F_2(z)&=\chi(z)E_1(z)=\chi(z)  \f{\left(A (z) + S_1\right)^{-1}
-(A_0 + S_1)^{-1}}{z} \\
&=  \chi(z) \sum_{k=1}^\infty (-1)^k z^{k-1} (A_0 + S_1)^{-1}
\left[A_1 (z) (A_0+S_1)^{-1}\right]^k.
\end{align*}
Let $\chi_1$ be a smooth cut off function which is equal to $1$ in $[-1,1]$. Note that  the support of
$\chi$ is contained in $[-1,1]$. We have
$$
F_2(z)=  \sum_{k=1}^\infty (-1)^k   \chi(z)z^{k-1} (A_0 + S_1)^{-1}
\left[  \chi_1(z)A_1(z) (A_0+S_1)^{-1}\right]^k.
$$
Using Lemma~\ref{L:convolve} and Young's inequality, we obtain
\begin{align} \label{F_2}
& \int_{-\infty}^\infty \left\| \left|\widehat{F_2}(\xi)\right|\right\|_{ 2\rightarrow 2}d\xi \leq \\
& \sum_{k=1}^\infty \|\widehat{(\chi(z)z^{k-1})}\|_{L^1} \| |(A_0+S_1)^{-1}| \|_{2\rightarrow 2}^{k+1}
\left[\int_{-\infty}^\infty \| | \widehat{(\chi_1 A_1 )}(\xi)| \|_{2\rightarrow 2} d\xi\right]^k. \nonumber
\end{align}
By an argument similar to Remark~1 in \cite{ErdSch} it is easy to see that
 $|(A_0+S_1)^{-1}|$ is bounded on $L^2$. Also note that
\begin{align}
\|\widehat{(\chi(z)z^{k-1})}\|_{L^1}& \lesssim \|(1+|\xi|)\widehat{(\chi(z)z^{k-1})}(\xi)\|_{L^2}\nonumber\\
&\lesssim \|\chi(z)z^{k-1}\|_2 +\| \f{d}{dz}(\chi(z)z^{k-1}) \|_2
 \lesssim \lambda_0^k. \label{chilambda}
\end{align}
Below, we prove that
\bea\label{chia1}
\int_{-\infty}^\infty \| | \widehat{(\chi_1 A_1 )}(\xi)| \|_{2\rightarrow 2} d\xi\lesssim 1.
\eea
If $\lambda_0$ is chosen sufficiently small, using  \pref{chilambda} and \pref{chia1} in \pref{F_2} completes the proof of the
lemma for $F_2$.
Recall that
\begin{align*}
A_1(z)(x,y) &=
\f{1}{4\pi z|x-y| }
v_2(x)\left[ \begin{array}{cc}
 e^{iz|x-y|}-1  & 0 \\
 0& -e^{-\sqrt{2\mu+z^2}|x-y|} +e^{-\sqrt{2\mu}|x-y|} \end{array} \right] \, v_1(y)\\
&= \f{1}{4\pi  }v_2(x)
\left[ \begin{array}{cc}
 i \int_0^1 e^{iz|x-y|b} db  & 0 \\
 0&  \int_0^1 \f{bz}{\sqrt{2\mu+z^2b^2}} e^{-\sqrt{2\mu+z^2b^2}|x-y|}  db \end{array} \right]
    \,v_1(y).
\end{align*}
We have (with $h_a(z)=\f{ z}{\sqrt{2\mu a^2+z^2 }} e^{-\sqrt{2\mu a^2+z^2 } }$)
\begin{align*}
& \widehat{(\chi_1 A_1 )}(\xi)(x,y)  =
\\
& \f{v_2(x)}{4\pi}
\left[\begin{array}{cc}
 i \int_0^1 \widehat{\chi_1}(\xi-|x-y|b) db  & 0 \\
 0&  \int_0^1 \int_{-\infty}^\infty \widehat{\chi_1}(\xi-b |x-y|\eta ) \widehat{h_{|x-y|}}(\eta) d\eta db   \end{array} \right]
  v_1(y).
\end{align*}
Hence by Schur's test, we can bound
$
\int_{-\infty}^\infty  \|  |  \widehat{\chi_1 A_1}(\xi)   | \|_{2\rightarrow 2} \; d\xi$ by a sum of quantities
of the form
\begin{align} \label{bitartik}
& \int_{-\infty}^\infty \sup_x  \int_{\R^3} \int_0^1     | \widehat{\chi_1}(\xi-|x-y|b)|
 |w_1(x)||w_2(y)|  db \, dy  \,d\xi,\, \text{ and }\\
&  \int_{-\infty}^\infty \sup_x  \int_{\R^3} \int_0^1    \int_{-\infty}^\infty | \widehat{\chi_1}(\xi-\eta|x-y|b)|
|\widehat{h_{|x-y|}}(\eta)|
 |w_1(x)||w_2(y)|  d\eta\, db \, dy  \,d\xi, \label{bitartikk}
\end{align}
where $w_1$ and $w_2$ satisfy
\bea \label{vxvy}
|w_1(x)||w_2(y)|\lesssim  \la x\ra^{-\beta/2} \la y\ra^{-\beta/2}\lesssim \la x-y \ra^{-\beta/2}.
\eea
Using (\ref{vxvy}) in (\ref{bitartik}),
we obtain
\begin{align*}
\pref{bitartik} &\lesssim \int_{-\infty}^\infty \sup_x  \int_0^1 \int_{\R^3} \la x-y \ra^{ -\beta/2}
 | \widehat{\chi_1}(\xi-|x-y|b)|
      dy  \, db\, d\xi \\
 & =  \int_{-\infty}^\infty \int_0^1 \int_{\R^3} \la  y \ra^{ -\beta/2}
 | \widehat{\chi_1}(\xi-| y|b)|
      dy  \, db\, d\xi \\
 &\leq \|\widehat{\chi_1}\|_1 \int_{\R^3} \la  y \ra^{ -\beta/2} dy  <\infty
\end{align*}
provided   $\beta>6$, i.e. $|V(x)|\lesssim \langle x \rangle^{-6-}$.
Now, we bound (\ref{bitartikk}).
Using (\ref{vxvy}) (with $\beta =0$)  and (\ref{hb}) in (\ref{bitartikk}),
we obtain
\begin{align*}
\pref{bitartikk} &\lesssim \int_{-\infty}^\infty \sup_x  \int_{-\infty}^\infty \int _0^1 \int_{\R^3} |\widehat{\chi_1}(\xi-\eta|x-y|b)|
     \f{e^{-\sqrt{\mu}|x-y|}}{|x-y|\la \eta \ra^2}   dy  \, db \,  d\eta \,d \xi \\
 & =  \int_{-\infty}^\infty  \int_{-\infty}^\infty \int _0^1 \int_{\R^3} |\widehat{\chi_1}(\xi-\eta|y|b)|
     \f{e^{-\sqrt{\mu}| y|}}{| y|\la \eta \ra^2}   dy  \, db \,  d\eta \,d \xi\\
&\leq  \|\widehat{\chi_1}\|_1  \int_{-\infty}^\infty   \int_{\R^3}
     \f{e^{-\sqrt{\mu}| y|}}{| y|\la \eta \ra^2}   dy  \,    d\eta
  < \infty.
\end{align*}

Next, we consider  $F_4$:
$$
F_4(z) = \chi(z) S_1 E_2(z) S_1
=\chi(z) S_1 \sum_{k=1}^\infty (-1)^k z^{k-1}  m (0)^{-1}
\left[m_1(z) m (0)^{-1}\right]^k S_1.
$$
Arguing as in the case of $F_2$, it suffices to prove that \beeq
\label{eq:m1FT} \int_{-\infty}^\infty \| | \widehat{(\chi_1 m_1
)}(\xi)| \|_{2\rightarrow 2} d\xi\lesssim 1, \eneq where $\chi_1$
is a smooth cut-off function which is equal to $1$ in the support
of $\chi$ (i.e. in $[-\lambda_0,\lambda_0]$) and which is
supported in $[-\lambda_1,\lambda_1]$. Recall that
$$
m_1(z)=S_1 \f{A_1(z)- A_1 (0)}{z} S_1
+\sum_{j=1}^\infty S_1 (-1)^j z^{j-1}
\left(A_1(z) (A_0 + S_1)^{-1}\right)^{j+1} S_1.
$$
The second summand can be analyzed as above
(here $\lambda_1$ is chosen sufficiently small to guarantee the convergence of the series, and than we choose $\lambda_0$ even smaller).
Now, we consider the first summand. Note that
\begin{align} \label{eq:A2def}
& A_2(z)(x,y) := \f{A_1(z) - A_1(0)}{z}(x,y)   \\
&= \f{1}{4\pi  }v_2(x)
\left[ \begin{array}{cc}
 i \int_0^1 \f{e^{iz|x-y|b}-1}{z} db  & 0 \\
 0&  \int_0^1 \f{b}{\sqrt{2\mu+z^2b^2}} e^{-\sqrt{2\mu+z^2b^2}|x-y|}  db \end{array} \right]
    \,v_1(y)  \nn
\\
    &= \f{1}{4\pi  }v_2(x)
\left[ \begin{array}{cc}
 - \int_0^1  |x-y| (1-b) e^{iz|x-y|b}  db  & 0 \\
 0&  \int_0^1 \f{b}{\sqrt{2\mu+z^2b^2}} e^{-\sqrt{2\mu+z^2b^2}|x-y|}  db \end{array} \right]
    \,v_1(y)  \nn
\end{align}
This can be analyzed as in the previous case. Because of the additional $|x-y|$ term, we need  to have
 $\beta>8$, i.e. $|V(x)|  \lesssim \langle x \rangle ^{-8-}$.

Next, we  deal with   $\f{d}{dz}F_j(z)$. Once again we omit the analysis of $F_1$ and $F_3$.
Note that
\begin{align*}
\f{d}{dz}F_2(z) = &   \sum_{k=1}^\infty (-1)^k \f{d}{dz}\left(\chi(z) z^{k-1}\right) (A_0 + S_1)^{-1}
\left[A_1(z) (A_0+S_1)^{-1}\right]^k \\
 +& \sum_{k=1}^\infty (-1)^k \chi(z) z^{k-1} (A_0 + S_1)^{-1}\times\\
&\times\sum_{j=1}^k  [A_1(z) (A_0+S_1)^{-1} ]^{j-1}  [\f{d}{dz}A_1(z) (A_0+S_1)^{-1}]
 [A_1(z) (A_0+S_1)^{-1}]^{k-j}
\end{align*}
Arguing as above, it suffices to prove that
\bea\label{chia1prime}
\int_{-\infty}^\infty \| | \widehat{(\chi_1(A_1)^\prime)}(\xi)| \|_{2\rightarrow 2} d\xi\lesssim 1.
\eea
Note that
(with $a=|x-y|$)
\begin{align*}
&\f{d}{dz}A_1(z) (x,y)
 =  \f{1}{4\pi  }v_2(x)  \times \left[\begin{array}{cc}
 -|x-y|\int_0^1 b e^{iz|x-y|b} db  & 0 \\
 0&  \int_0^1 \f{d}{dz} h_a(abz) db \end{array} \right]
    \,v_1(y)
\end{align*}
These are similar to the terms treated above. Therefore \pref{chia1prime} holds
provided $|V(x)|  \lesssim \langle x \rangle ^{-8-}$.

Finally, we analyze $\f{d}{dz}F_4(z)$. In view of the preceding, it suffices to prove that
\bea\label{a2prime}
\int_{-\infty}^\infty \| | \widehat{(\chi_1(A_2)^\prime)}(\xi)| \|_{2\rightarrow 2} d\xi\lesssim 1.
\eea
We have
\begin{align*}
&\f{d}{dz}A_2(z)(x,y) =
\\
& \f{1}{4\pi  }v_2(x)
\left[ \begin{array}{cc}
 -i \int_0^1  |x-y|^2 (1-b)b e^{iz|x-y|b}  db  & 0 \\
 0&  \int_0^1 \f{d}{dz}\left[\f{b}{\sqrt{2\mu+z^2b^2}} e^{-\sqrt{2\mu+z^2b^2}|x-y|}\right]  db \end{array} \right]
    \,v_1(y)  \nn
\end{align*}
These are treated as before; \pref{a2prime} holds provided $|V(x)|  \lesssim \langle x \rangle ^{-10-}$.
\end{proof}

Now, we consider $K_3^{12}$. We omit the analysis of the other components of $K_3$ since they can be handled similarly.
Denote
$$\int_{\R^6}\f{d}{dz}\left(\chi(z)e^{iz |x-u_1|-\sqrt{2\mu+z^2}|y-u_2| }
(I_1v_1E(z)v_2I_2^T)(u_1,u_2)\right) \f{d u_1 d u_2}{|x-u_1||y-u_2|}  $$
by $\calG_{x,y}(z)$.
By Lemma~\ref{L:birb}, it suffices to prove (\ref{fl1}) for $\calG_{x,y}$.
Let us concentrate on the term where the derivative hits $\chi (I_1v_1E v_2I_2^T)$ (the term
where the derivative hits the exponential is similar):
$$
\tilde{\calG}_{x,y}(z)=\int_{\R^6}[\chi(z)(I_1v_1E(z)v_2I_2^T)]^\prime(u_1,u_2)
\f{e^{iz |x-u_1|-\sqrt{2\mu+z^2}|y-u_2| }}{|x-u_1||y-u_2|}du_1du_2.
$$
Similarly (we denote $e^{-a\sqrt{2\mu+z^2}}$ by $e_a$)
\begin{align}
\|\widehat{\tilde{\calG}_{x,y}}(\xi)\|_{L^1}&\leq\int_{\R^8}
\left|\widehat{ [\chi(I_1v_1Ev_2I_2^T)]^\prime }(\eta)(u_1,u_2)\widehat{e_{|y-u_2|}}(\xi-\eta-|x-u_1|)
\right|
\f{ d\eta
d\xi du_1du_2 }{|x-u_1||y-u_2|}\nn\\
&\leq   \sup_a \|\widehat{e_{a}}\|_1  \int_{\R^7} \left|\widehat{ [\chi(I_1v_1Ev_2I_2^T)]^\prime }(\eta)(u_1,u_2)\right|
\f{ d\eta du_1du_2 }{|x-u_1||y-u_2|} \nn\\
 & \lesssim \sup_a \|\widehat{e_{a}}\|_1  \sum_{i,j=1}^2
\int_{-\infty}^\infty \left\|\left|\widehat{(\chi E_{ij} )^\prime}(\eta) \right| \right\|_{L^2\rightarrow L^2}d\eta.
\nn
\end{align}
It is not difficult to see that (using Lemma~\ref{lem:calc}) $ \sup_a \|\widehat{e_{a}}\|_1 <\infty$.
Therefore, for $\tilde{\calG}_{x,y}$, (\ref{fl1}) follows from
(\ref{eq:K3end}).

\subsection{The general case}

We now prove Theorem~\ref{T:scalar2} in the general case.
Using (\ref{eq:teilung})  in (\ref{eq:new}), we have
\begin{align} \label{eq:MR_S_2}
\mr(z)& = \mr_0(z) - \mr_0(z) v_1 \Gamma_1(z) v_2 \mr_0(z) \\
 & \quad -\frac{1}{z}\mr_0(z) v_1 \Gamma_1(z)S_1\Gamma_2(z)S_1\Gamma_1(z) v_2
 \mr_0(z) \nn\\
 & \quad - \frac{1}{z^2} (\mr_0(z) T(z) \mr_0(z) - \mr_0(0) T(0) \mr_0(0))\nn\\
  & \quad - \frac{1}{z^2} \mr_0(0) T(0) \mr_0(0),\nn
\end{align}
where $\Gamma_1(z)=(A(z)+S_1)^{-1}$,
$\Gamma_2(z)=(m(z)+S_2)^{-1}$ and
\[
T(z)= v_1 \Gamma_1(z)S_1 \Gamma_2(z) S_2 b(z)^{-1}
S_2 \Gamma_2(z) S_1 \Gamma_1(z) v_2.
\]
Substituting (\ref{eq:MR_S_2}) in (\ref{Kt}), we have (ignoring $2e^{it\mu}$)
\begin{align}
K_t &= \int_{-\infty}^\infty e^{it z^2}  \chi(z) z  \mr_0(z)dz -
\int_{-\infty}^\infty e^{it z^2}  \chi(z) z \mr_0(z) v_1 \Gamma_1(z) v_2 \mr_0(z)  dz \label{ktS2_1}\\
& \quad -\int_{-\infty}^\infty e^{it z^2}\chi(z)  \mr_0(z) v_1 \Gamma_1(z)S_1\Gamma_2(z)S_1\Gamma_1(z) v_2 \mr_0(z)  dz \label{ktS2_2} \\
& \quad -\int_{-\infty}^\infty e^{it z^2}\chi(z) \frac{1}{z} (\mr_0(z) T(z) \mr_0(z) - \mr_0(0) T(0) \mr_0(0))   dz.  \label{ktS2_3}
\end{align}
Here, the singular term $\frac{1}{z^2} \mr_0(0) T(0) \mr_0(0)$ in (\ref{eq:MR_S_2}) has no contribution since
the integral in (\ref{Kt}) is a principal value integral and the integrand is odd.
The first operator in (\ref{ktS2_1}) is dispersive since it is the low energy part of the free evolution.
The second operator in (\ref{ktS2_1}) is also dispersive, which can be proved by repeating the analysis of $K_3$
in the previous section. The operator in (\ref{ktS2_2}) can be rewritten as a sum of two operators one similar
two $K_1$ and the other similar to $K_3$ in the previous section.
The $L^1 \rightarrow L^\infty$ norm of the former decays like $t^{-1/2}$ and the latter is dispersive.
Now, we consider (\ref{ktS2_3}). We can write it as a sum of the following operators:
\begin{align}
& \int_{-\infty}^\infty e^{it z^2}\chi(z)  \mr_0(z) \frac{T(z)-T(0)}{z} \mr_0(z)    dz,  \label{ktS2_31}\\
& \int_{-\infty}^\infty e^{it z^2}\chi(z) \frac{\mr_0(z)-\mr_0(0)}{z}  T(0) \mr_0(z)  dz  \label{ktS2_32},\\
& \int_{-\infty}^\infty e^{it z^2}\chi(z)  \mr_0(0) T(0) \frac{\mr_0(z) - \mr_0(0)}{z}    dz  \label{ktS2_33}.
\end{align}
Since we don't have an extra power of $z$, the $L^1 \rightarrow L^\infty$ norm of
these operators decay like  $t^{-1/2}$ (see Lemma~\ref{L:birb}).
First let us consider (\ref{ktS2_32}). Note that
\[
T(0)(x,y)=v_1(x) (S_2b(0)^{-1}S_2)(x,y)v_2(y)
\]
is a finite rank operator. Therefore it suffices to study operators with kernel (with the notation $a_1=|x-u_1|$, $a_2=|y-u_2|$)

\begin{align*}
 &  \int_{\R^7}  \frac{e^{it z^2}\chi(z)}{a_1a_2}
\left[
\begin{array}{ll}
\frac{e^{iza_1}-1}{z} & 0\\0& \frac{e^{-\sqrt{2\mu}a_1}-e^{-\sqrt{2\mu+z^2}a_1}}{z}
\end{array}
\right]   \\
& \quad
\left[
\begin{array}{l}
r_1(u_1)  \\ r_2(u_1)
\end{array}
\right]
\left[
\begin{array}{l}
c_1(u_2)  \\ c_2(u_2)
\end{array}
\right]^T
\left[
\begin{array}{ll}
 e^{iza_2}  & 0\\0&  -e^{-\sqrt{2\mu+z^2}a_2}
\end{array}
\right]
du_1du_2dz,
\end{align*}
where $r_1, r_2, c_1, c_2 \in  L^{2,\beta/2}$.
This can be rewritten as (with the notation
$h_a(z)=\frac{z}{\sqrt{2\mu a^2+z^2}}e^{-\sqrt{2\mu a^2+z^2}}$)
\begin{align*}
 & \int_{\R^7} \int_0^1 \frac{e^{it z^2}\chi(z)}{ a_2}   \\
& \quad
\left[
\begin{array}{ll}
-i  r_1(u_1) c_1(u_2) e^{iz(ba_1+a_2)}  & i r_1(u_1) c_2(u_2) e^{ibza_1} e^{-\sqrt{2\mu+z^2}a_2}  \\
 r_2(u_1) c_1(u_2) h_{a_1}(a_1bz)e^{iza_2} &
-  r_2(u_1) c_2(u_2)h_{a_1}(a_1bz) e^{-\sqrt{2\mu+ z^2} a_2 }
\end{array}
\right]
dbdu_1du_2dz.
\end{align*}
This operator is similar to the operator $K_1$ studied in the previous section. We omit the analysis.
Now, we consider (\ref{ktS2_32}). Similarly it suffices to consider operators of the form:
\begin{align*}
  &  \int_{\R^7}  \frac{e^{it z^2}\chi(z)}{a_1a_2}
\left[
\begin{array}{ll}
1 & 0\\0&  -e^{-\sqrt{2\mu}a_1}
\end{array}
\right]   \\
&\quad
\left[
\begin{array}{l}
r_1(u_1)  \\ r_2(u_1)
\end{array}
\right]
\left[
\begin{array}{l}
c_1(u_2)  \\ c_2(u_2)
\end{array}
\right]^T
  \left[
  \begin{array}{ll}
   \frac{e^{iza_2}-1}{z}  & 0\\0&  \frac{e^{-\sqrt{2\mu }a_2}-e^{-\sqrt{2\mu+z^2}a_2}}{z}
  \end{array}
   \right]    du_1du_2dz \\
&=\int_{\R^7} \int_0^1 \frac{e^{it z^2}\chi(z)}{a_1 }
\left[
\begin{array}{ll}
  -i r_1(u_1) c_1(u_2)  e^{ibza_2} &  - r_1(u_1) c_2(u_2) h_{a_2}(a_2bz)     \\
i r_2(u_1) c_1(u_2) e^{-\sqrt{2\mu}a_1} e^{ibza_2}  &
 r_2(u_1) c_2(u_2)e^{-\sqrt{2\mu}a_1} h_{a_2}(a_2bz)
\end{array}
\right]
dbdu_1du_2dz.
\end{align*}
Once again this operator is similar to $K_1$ studied in the previous section.
Now, we consider (\ref{ktS2_31}). We use the following identity
\begin{align*}
 T(z)-T(0) &= v_1 (\Gamma_1(z)-\Gamma_1(0)) S_1\Gamma_2(z)S_2b(z)^{-1}S_2\Gamma_2(z)S_1\Gamma_1(z)   v_2 \\
&+ v_1  S_1(\Gamma_2(z)-\Gamma_2(0)) S_2b(z)^{-1}S_2\Gamma_2(z)S_1\Gamma_1(z)   v_2 \\
&+ v_1   S_2(b(z)^{-1}-b(0)^{-1})S_2\Gamma_2(z)S_1\Gamma_1(z)   v_2  \\
&+ v_1   S_2 b(0)^{-1} S_2(\Gamma_2(z)-\Gamma_2(0))S_1\Gamma_1(z)   v_2 \\
&+ v_1   S_2 b(0)^{-1} S_2 (\Gamma_1(z)-\Gamma_1(0)) v_2.
\end{align*}
In view of the analysis of the operator $K_3$ in the previous section, it suffices to prove the
  bound  (\ref{hammal}) for the following basic building blocks:
\begin{align*}
 F_1(z) &= \chi(z)\Gamma_1(z)=\chi(z)(A(z)+S_1)^{-1} \\
 F_2(z) &= \chi(z)z^{-1}(\Gamma_1(z)-\Gamma_1(0))=\chi(z)z^{-1}((A(z)+S_1)^{-1}-(A_0+S_1)^{-1}) \\
 F_3(z) &= \chi(z)S_1\Gamma_2(z)S_1=\chi(z)S_1(m(z)+S_2)^{-1}S_1 \\
 F_4(z) &= \chi(z)z^{-1}S_1(\Gamma_2(z)-\Gamma_2(0))S_1
 =\chi(z)z^{-1}S_1((m(z)+S_2)^{-1}-(m(0)+S_2)^{-1})S_1
\\
F_5(z) &= \chi(z) S_2 b(z)^{-1} S_2 = \chi(z) S_2 (b(0)+z b_1(z))^{-1} S_2 \\
F_6(z) &= \chi(z) S_2 z^{-1}(b(z)^{-1}-b(0)^{-1}) S_2.
\end{align*}
The functions $F_j$ with $1\le j\le 4$ were already discussed in
Lemma~\ref{L:hammal}. Therefore, it suffices to prove that
\beeq
\label{eq:F56} \max_{j=5,6}\int_{-\infty}^\infty \left\|\,
\left|\widehat{F_j }(\xi)\right|\, \right\|_{ 2\rightarrow 2}d\xi
< \infty.
\eneq

Recall that, see (\ref{eq:bdef}),
\begin{align}
 b(0) &= S_2 m_1(0) S_2 \nn \\
 b(z) &= b(0) + z b_1(z) = b(0)(1+z b(0)^{-1}b_1(z)) \nn \\
 b_1(z) &=  \f{S_2 [m_1(z)-m_1(0)]S_2}{z} +
\f{1}{z} \sum_{k=1}^\infty (-1)^k z^k S_2\left(m_1(z) (m(0)+S_2)^{-1}\right)^{k+1}S_2 \label{eq:m1ser}\\
 b(z)^{-1} &= \sum_{j=0}^\infty (-1)^j z^j (b(0)^{-1}b_1(z))^j b(0)^{-1}. \label{eq:b1ser}
\end{align}
Applying Lemma~\ref{L:convolve}  to the Neuman series in~\eqref{eq:b1ser} shows that in order to obtain~\eqref{eq:F56},
we need to prove that
\[
\int_{-\infty}^\infty \left\|\, \left|\widehat{\chi_1 b_1 }(\xi)\right|\,\right\|_{ 2\rightarrow 2}d\xi < \infty.
\]
Another application of Lemma~\ref{L:convolve}, this time to the
Neuman series~\eqref{eq:m1ser}, reduces matters to proving
\[
\int_{-\infty}^\infty \left\|\, \left|\widehat{\chi_2 m_1 }(\xi)\right|\,\right\|_{ 2\rightarrow 2}d\xi < \infty,
\]
which was already done in \eqref{eq:m1FT}. In both these cases,
the cut-off functions $\chi_1, \chi_2$ need to be taken with
sufficiently small supports. This leaves the term
\[
\f{S_2 [m_1(z)-m_1(0)]S_2}{z}
\]
from~\eqref{eq:m1ser} to be considered.
In view of
\eqref{eq:m1def} and \eqref{eq:A2def},
\begin{align*}
& S_2\frac{m_1 (z)-m_1(0)}{z}S_2 \\
& = S_2
\f{A_2 (z) - A_2 (0)}{z}S_2
+\sum_{k=1}^{\infty} (-1)^k z^{k-1 } S_2
\left(A_1(z) (A_0 + S_1)^{-1}\right)^{k+1}S_2.
\end{align*}
By \eqref{chia1}, and Lemma~\ref{L:convolve}, the Neuman series
makes a summable contribution to~\eqref{eq:F56}. On the other
hand, the contribution of
\[
S_2\f{A_2 (z) - A_2 (0)}{z}S_2
\]
to \eqref{eq:F56} is controlled by the bound~\eqref{a2prime}, and we are done.


\end{document}